\documentclass[12pt, russian]{article}
\usepackage[T2A]{fontenc}
\usepackage[cp866]{inputenc}
\usepackage[english]{babel}
\usepackage{amsmath,latexsym, amsfonts, amssymb, amscd,amsthm}
\date{}
\DeclareMathOperator{\di}{div}
\begin{document}
\title{On a new type of the $\ell$-adic regulator for algebraic number fields
(The $\ell$-adic regulator without logarithms)}
\author{L.~V.~Kuz'min\thanks{This work was done with the financial support of the Russian Foundation for basic Reseach (grant no. 11-01-00588-a)}} 
\maketitle
\newtheorem{prop}{Proposition}[section]
\newtheorem{theorem}{Theorem}
\newtheorem{lemma}{Lemma}[section]
\newtheorem{conj}{Conjecture}
\newtheorem*{con}{Consequence}
\newtheorem*{co}{The main conjecture}
\newtheorem*{de}{Definition}
\newcommand{\lr}{\longrightarrow}
\renewcommand{\theequation}{\thesection.\arabic{equation}}
\begin{abstract}
For an algebraic number field $K$ such that prime $\ell$ splits completely in $K$
we define a regulator $\mathfrak R_\ell(K)\in\mathbb Z_\ell$ that characterize the subgroup of
universal norms from the cyclotomic $\mathbb Z_\ell$-extension of $K$ in the
completed group of $S$-units of $K$, where $S$ consists of all prime divisors of
$\ell$. We prove that inequality $\mathfrak R_\ell(K)\neq0$ follows from the
$\ell$-adic Schanuel conjecture and holds true for some Abelian extensions 
of imaginary quadratic fields.
\end{abstract}
\section{Introduction}
Let $K$ be an algebraic number field and $\ell$ a fixed prime number. Then there
is the $\ell$-adic regulator $R_\ell(K)\in \mathbb Q_\ell$ of $K$ in the sense
of~\cite{reg1}, which characterize the  $\ell$-adic behavior of the group of
units $U(K)$ of $K$. (Its definition is reproduced in Section 2, see Conjecture
3.) Nevertheless this regulator isn't perfect from the viewpoint of Iwasawa
theory. 

For example, suppose that prime $\ell$ splits completely in $K$. Let
$K_\infty$ be the cyclotomic $\mathbb Z_\ell$-extension of $K$. In this case for
any intermediate subfield $K_n$ of the extension $K_\infty/K$, where
$[K_n:K]=\ell^n$ and $K=K_0$, one has a regulator $R_\ell(K_n)$, which
characterize the group of units $U(K_n)$ of $K_n$.  As it was shown in
\cite{reg5}, the behavior of the $\ell$-adic exponents $\nu_\ell(R_\ell(K_n))$
for increasing $n$ shows some interesting feathers. This event is connected with
the fact that the properties of the sequence of the groups of units
$\{U(K_n)\}_{n\geqslant0}$ reflect in a single object -- a Galois module (with
respect to the action of the Galois group $\Gamma:= G(K_\infty/K)$) $\bar
U\bigl(K_\infty):=\varprojlim (U(K_n)[\ell]/\mu_\ell(K_n)\bigr)$, where
$U(K_n)[\ell]$ is the pro-$\ell$-completion of $U(K_n)$, $\mu_\ell(K_n)$ is the
group of all roots of unity of $\ell$-power degree in $K_n$ and the projective
limit is taken with respect to the norm maps. In this case for any $n$ one has a
natural projection $\pi_n:\bar U(K_\infty)\to U(K_n)[\ell]/\mu_\ell(K_n)$. In
particular, one has a projection $\pi:\bar U(K_\infty)\to
U(K)[\ell]/\mu_\ell(K)$. But if prime $\ell$ splits completely in $K$ and
Leopoldt conjecture holds for $K$ and $\ell$, then one has $\pi=0$ (we should
explain it more carefully in Section 2).
So it is natural to consider instead of $U(K)[\ell]/\mu_\ell(K)$ some other
object, which behaves better from the viewpoint of Iwasawa theory.

Namely, let $S$ be the set of all places over $\ell$ in any field under
consideration and $U_S(K_n)$ the group of $S$-units in $K_n$. Then one can
link the sequence of groups $\{U_S(K_n)\}_{n\geqslant 0}$ to a single
$\Gamma$-module $\bar
U_S(K_\infty):=\varprojlim\bigl(U_S(K_n)[\ell]/\mu_\ell(K_n)\bigr)$, where the
limit is taken with respect to the norm maps. For any $n\geqslant 0$, there is a
natural projection $\pi_n':\bar U_S(K_\infty)\to U_S(K_n)[\ell]/\mu_\ell(K_n)$.
The compactness of the groups $U_S(K_n)[\ell]/\mu_\ell(K_n)$ yields that the
image of $\pi'_n$, which we denote by $U_{S,1}(K_n)$, coincides with the
subgroup of universal norms in $U_S(K_n)[\ell]/\mu_\ell(K_n)$, that is,
\begin{equation}
 U_{S,1}(K_n)=\bigcap_{m>n}N_{K_m/K_n}\bigl(U_S(K_m)[\ell]/\mu_\ell(K_m)\bigr).
\label{1.0}
\end{equation}
The new $\ell$-adic regulator, which we wish to defined, should characterize the
$\ell$-adic behavior of the elements of $U_{S,1}(K_n)$. But in the present paper
we shall treat only the case $n=0$. Thus we shall define the new $\ell$-adic
regulator $\mathfrak R_\ell(K)\in \mathbb Z_\ell$ that characterize $U_{S,1}(K)$
under assumption that prime $\ell$ splits completely in $K$. (The case of
arbitrary $K$ will be treated in the next paper.) This restriction enable us to
make the presentation more transparent.

Let $D=D(K)$ be  a free additive $\mathbb Z_\ell$-module generated by the
elements of $S$. Then one can put to any $a\in U_S(K)[\ell]/\mu_\ell(K)$ its
``principal divisor'' $\di (a)\in D$. Thus one gets a homomorphism
$\di:U_{S,1}(K)\to D$, which under assumption of validity of Leopoldt conjecture
for $K$ and $\ell$, induces an injection $\di :U_{S,1}(K)\hookrightarrow D$.

There is a standard $\mathbb Z_\ell$-linear scalar product on $D$
\begin{equation}
\langle\phantom x,\phantom x\rangle:D\times D\to\mathbb Z_\ell.
\label{1.1}
\end{equation}
Namely, if $S=\{v_1,\ldots,v_m\}$ and $\alpha,\beta\in
D,\quad\alpha=\sum_{i=1}^ma_iv_i,\quad b=\sum_{i=1}^m b_iv_i$, where $a_i,b_i\in
\mathbb Z_\ell$, then we put $\langle\alpha,\beta\rangle=\sum_{i=1}^ma_ib_i$.
Assuming the Leopoldt conjecture valid for $K$ and $\ell$, we prove that
$U_{S,1}(K)$ is a free $\mathbb Z_\ell$-module of rank $t=r_1+r_2$, where $r_1$
and $r_2$ is the number of real and complex places of $K$ respectively. Then we
define a new regulator $\mathfrak R_\ell(K)\in \mathbb Z_\ell$ by the formula
\[
 \mathfrak R_\ell(K)=\det(\langle \di(e_i),\di(e_j)\rangle),\quad 1\leqslant
i,j\leqslant t,
\]
where $e_1,\ldots,e_t$ is a $\mathbb Z_\ell$-basis of $U_{S,1}(K)$. 

Note that $\mathfrak R_\ell(K)$ is defined only up to multiplication by the
square of an arbitrary element of $\mathbb Z_\ell^\times$.

In the present paper we consider the following
\begin{co} Let $K$ be an algebraic number field such that prime $\ell$ splits
completely in $K$. Then $\mathfrak R_\ell(K)\neq 0$.
\end{co}
Obviously, the main conjecture is equivalent to assertion that the product
\eqref{1.1} induces a non-degenerate scalar product on $U_{S,1}(K)$, which is
defined by
\begin{equation}
 \{\alpha,\beta\}=\langle\di(\alpha),\di(\beta)\rangle
\label{1.2}
\end{equation}
for any $\alpha,\beta\in U_{S,1}(K)$.

We prove that the main conjecture follows from the $\ell$-adic Schanuel
conjecture (see Conjecture 1 in Section~2). Then we give an unconditional proof
of the main conjecture for some types of Abelian extensions of imaginary quadratic
fields.

Now we shall give the plan of the paper and describe our main results in more
details.

In Section~2 we remind some necessary definitions and formulate the conjectures that
we will use in sequel. Some of them, such as Conjectures 1 and 2) are
well-known, while Conjectures 3 and 4 were firstly formulated by the author in
\cite{reg1}. The main result of Section~2 is Theorem~1, which gives the $\mathbb
Z_\ell$-rank of $U_{S,1}(K)$. Moreover, if $K/\mathbb Q$ is a Galois extension
with Galois group $G$ and such that the Leopoldt conjecture holds for a pair $(K,\ell)$,
then Theorem~1 asserts that the modules $\bigl(U(K)[\ell]\otimes_{\mathbb
Z_\ell}\mathbb Q_\ell\bigr)\oplus\mathbb Q_\ell$ and $U_{S,1}(K)\otimes_{\mathbb
Z_\ell}\mathbb Q_\ell$ are isomorphic over $\mathbb Q_\ell[G]$. (Here and in
many cases later we write the multiplication additively).

If $K/\mathbb Q$ is not Galois, then Theorem~1 asserts only the equality of
dimensions of these two $\mathbb Q_\ell$-spaces. The proof of this theorem
depends heavily on class field theory and Iwasawa theory. We shall use also an
important consequence of Theorem 1 on an existence of a special element in the
group $U_{S,1}(K)$, which is analogous to the Artin unit in the group of units
$U(K)$. Note that we don't assume that $\ell$ splits completely in $K$.

In Section~3 we use this consequence of Theorem 1 to produce a conditional proof
of the main conjecture.  Namely, assuming that Conjecture~1 holds, we prove that
$\mathfrak R_\ell(K)\neq0$ for any algebraic number field $K$ and any $\ell$
that splits completely in $K$. In this proof one can assume that
$K/\mathbb Q$ is Galois with Galois group $G$ and $K$ contains an imaginary
quadratic field $k$. To do this, we introduce another conjecture (Conjecture 5),
which is similar to Leopoldt conjecture. But in distinction to the latter, it is
formulated in terms of a Galois module generated by $\alpha\in K$ such that its
principal divisor $(\alpha)$ equals to $\mathfrak l^{\bf h}$, where $\mathfrak
l$ is a prime divisor of $\ell$ in $K$ and ${\bf h}$ is a natural number. (The
class number of $K$, for example.) This conjecture is a consequence of
Conjecture~1, but we hope that it may be more tractable in some cases.

Using Conjecture 5, maп can define explicitly some group  $U_{S,2}(K)$, which
contains $U_{S,1}(K)$. The group $U_{S,2}(K)$ is the group of local universal
norms in $U_S(K)[\ell]/\mu_\ell(K)$. Using Conjecture 1, man can prove that the
index $(U_{S,2}(K):U_{S,1}(K))$ is finite (Proposition 3.5). It should be
mentioned that the last result gives a conditional proof of Conjecture~2 in
\cite{mt} for the case, when $\ell$ splits completely in $K$.

During calculation of the group $U_{S,2}(K)$, there appears some system of
coefficients $\{c_h\}_{h\in H}$, where $H=G(K/k)$ (see \eqref{3.7a}), which
enter the relations between the $\ell$-adic logarithms of some Artin unit
$\varepsilon$ in $U(K)$ and its conjugate, and the $\ell$-adic logarithms of
$\sigma(\alpha),\quad \sigma\in G$. We prove modulo Conjecture~1 (Lemma 3.1),
that after omitting from the system $\{c_h\}$ any element $c_{h_0}$, the rest
elements form an algebraically independent set. Then we introduce an analog of
the regulator $\mathfrak R_\ell(K)$, which characterizes not $U_{S,1}(K)$, but
the group $U_{S,2}(K)$. We prove that this analog is a non-zero polynomial in
$\{c_h\}_{h\in H, h\neq h_0}$. It proves Theorem~2, according to which
$\mathfrak R_\ell(K)\neq0$, if Conjecture~1 holds. 

In Section~4 we give an unconditional proof of the main conjecture for one
particular case (Theorem~3). We assume that $K$ is Galois over $\mathbb Q$ and
Abelian over an imaginary quadratic field $k$. We assume also that the only
non-trivial automorphism $\tau$ of $k/\mathbb Q$ acts on the Galois group
$H=G(k/\mathbb Q)$ by inversion, i.~e. $\tau h\tau^{-1}=h^{-1}$ for any $h\in
H$. We assume also that any irreducible character $\chi$ of $H$ with values in
$\bar {\mathbb Q}_\ell$ is conjugated with  $\chi^{-1}$ over $\mathbb Q_\ell$.

The proof is based on the same idea as the proof of Theorem~3.1 in \cite{reg1}.
Namely, instead of the product \eqref{1.2} we consider some more complicated
product~\eqref{4.0}, which has some additional symmetry (Proposition~4.1). If
the latter product is non-generate, we get immediately that $\mathfrak
R_\ell(K)\neq0$, nevertheless we can state this only under additional assumption
that $\chi$ and $\chi^{-1}$ are conjugate for any $\chi\in \hat H$.

The last condition depends only on the pair of numbers $(\ell,m)$, where $m$ is
an exponent of $H$. We give a full list of all such pairs in Proposition~4.3. One
can put a question in another way. Suppose $K$ to be fixed. What is the set of
all primes $\ell$ such that Theorem~3 holds for $(K,\ell)$, that is, $\mathfrak
R_\ell(K)\neq0$? The answer is given in Proposition~4.4.

A part of this paper was written when the author was a guest of Hausdorff
research institute for mathematics, Bonn (Trimester Program ``Arithmetic and Geometry'' 2013). 
The author would like to thank the collaborators of the institute for hospitality and the participants of the Trimester for stimulating communion.

\section{Preliminary results}
\setcounter{equation}{0}
Let $\bar{\mathbb Q}_\ell$ be the algebraic closure of the rational $\ell$-adic
number field $\mathbb Q_\ell$ and $\log:\bar{\mathbb Q}_\ell^\times\to
\bar{\mathbb Q}_\ell$ the $\ell$-adic logarithm in the sense of Iwasawa. Note
that in this paper we use neither real nor complex logarithms. So the denotation
$\log$ for the $\ell$-adic logarithm would not't cause any ambiguity. Thus, $\log$
is a homomorphism of the multiplicative group $\bar {\mathbb Q}_\ell^\times$
into the additive group $\bar {\mathbb Q}_\ell$. This homomorphism is uniquely
defined by conditions $\log(\ell)=0$, $\log(\xi)=0$ for any root of unity $\xi$,
whereas $\log$ is defined on the group of principal units $U^{(1)}(\bar{\mathbb
Q}_\ell)$ by the usual power series.

A well-known $\ell$-adic Schanuel conjecture asserts that for any
$x_1,\ldots,x_n\in\bar{\mathbb Q}_\ell^\times$ such that $\ell,x_1,\ldots,x_n$
are multiplicatively independent, that is the equality 
$\ell^c\prod_{i=1}^nx_i^{c_i}=1$ does not take place for any non-zero vector
$(c,c_1,\ldots,c_n)\in\mathbb Z^{n+1}$, the transcendence degree of $\mathbb
Q(x_1,\dots,x_n,\log x_1,\ldots,\log x_n)$ over $\mathbb Q$ is at least $n$. In
the present paper we shall use only the following particular case of this
conjecture.
\begin{conj}
\label{conj1}
 Let $x_1,\ldots,x_n\in\bar{\mathbb Q}_\ell^\times$ be algebraic over $\mathbb
Q$ and $\ell,x_1,\ldots,x_n$ are multiplicatively independent. Then $\log
x_1,\ldots,\log x_n$ are algebraically independent over $\mathbb Q$.
\end{conj}
Let $K$ be an algebraic number field and $\ell$ a fixed prime number. Then there
is a natural injection $i:K\hookrightarrow K\times_\mathbb Q\mathbb
Q_\ell\cong\prod_{i=1}^mK_{v_i}$, where $K_{v_i}$ is the completion of $K$ with
respect to $v_i$ and $v_1,\ldots,v_m$ run all the places of $K$ over $\ell$
(here we does not assume that $\ell$ splits completely in $K$). 

Considering $K_{v_i}$ as subfields of $\bar{\mathbb Q}_\ell$, we obtain
injections $\sigma_1,\ldots,\sigma_m$ of $K$ into $\bar{\mathbb Q}_\ell$ such
that $\sigma_i(K)\subset K_{v_i}$ and $i(a)=
\{\sigma_1(a),\ldots,\sigma_m(a)\}\in\prod_{i=1}^mK_{v_i}$ for any $a\in K$. In
particular, to any $x\in K^\times$ one can consider a vector
\begin{equation}
 \log_\ell(x)=\{\log \sigma_1(x),\ldots,\log\sigma_m(x)\}.
\label{2.1}
\end{equation}
Thus, we get a homomorphism
\[
 \log_\ell:K^\times\to\prod_{i=1}^m K_{v_i}.
\]
If $K$ is a Galois extension of $\mathbb Q$ with the Galois group $G$, then the
last homomorphism is a $G$-homomorphism. If $\ell$ splits completely in $K$,
one has $m=[K:\mathbb Q]$ and $\prod_{i=1}^m K_{v_i}\cong \mathbb Q_\ell^m$.
\begin{conj}
 \label{conj2}
{(Leopoldt conjecture)}. Let $u_1,\ldots,u_r\in U(K)$ be multiplicatively
independent units of $K$. Then $i(u_1),\ldots,i(u_r)$ are multiplicatively
independent over $\mathbb Z_\ell$ in $\prod_{i=1}^m K_{v_i}^\times$. In other
words, the vectors $\log_\ell u_1,\ldots,\log_\ell u_r$ are lineary independent
over $\mathbb Q_\ell$.
\end{conj}
The injection $i:U(K)\hookrightarrow\prod_{i=1}^mU(K_{v_i})$, where $U(K_{v_i})$
is the group of units of $K_{v_i}$, induces a homomorphism
$i':U(K)[\ell]\to\prod_{i=1}^mU(K_{v_i})$, where $U(K)[\ell]$ is the
pro-$\ell$-completion of $U(K)$. Conjecture \ref{conj2} is equivalent to the
assertion that $i'$ is a monomorphism. In general case, the group
$E_\ell(K):=\ker i'$ is known as Leopoldt kernel. So Leopoldt conjecture asserts
that $E_\ell(K)=1$ for any $(K,\ell)$. It is well known that Leopoldt conjecture
follows from the conjecture~\ref{conj1}. Leopoldt conjecture is proved for
Abelian fields and Abelian extensions of imaginary quadratic fields.

For an algebraic number field $K$ and a prime number $\ell$, let $K_\infty$ be
the cyclotomic $\mathbb Z_\ell$-extension of $K,\quad
G(K_\infty/K):=\Gamma\cong\mathbb Z_\ell$ the Galois group of $K_\infty/K$ and
$\gamma$ a fixed topological generator of $\Gamma$. By $K_n$ we denote the
unique intermediate subfield of  $K_\infty/K$ such that $[K_n:K]=\ell^n$. Hence,
$K=K_0$. We shall need some results on the group $E_\ell(K)$ and the sequence of
groups $E_\ell(K_n)$. The next statement, which is in fact a theorem, is known
as feeble Leopoldt conjecture.
\begin{prop}
 (See, for example, \cite{mt}). The groups $E(K_n)$ stabilize since some $n$,
that is, there is an index $n$, which depends only on $K$ and $\ell$, such that
the natural injection $E_\ell(K_n)\hookrightarrow E_\ell(K_{n_1})$ is an
isomorphism for any $n_1>n$.
\label{prop2.1}
\end{prop}
Let $M$ be the maximal Abelian $\ell$-extension of $K_\infty$ unramified outside
$\ell$ and $X=G(M/K_\infty)$. Then $\Gamma$ acts by conjugation on $X$, and thus
$X$ becomes a finitely generated $\Lambda$-module, where $\Lambda:=\mathbb
Z_\ell[[\Gamma]]$ is a completed group ring of $\Gamma$.
\begin{prop}
 (\cite{mt}, Proposition~5.4). There exists a natural isomorphism $X^\Gamma\cong
E_\ell(K)$.
\label{prop2.2}
\end{prop}
The next conjecture is a strengthening of Leopoldt conjecture. It is a
consequence of Conjecture~1 (see \cite[Appendix]{reg1}). We formulate it, since
it motivate many of our later considerations.
\begin{conj} 
(Conjecture on the $\ell$-adic regulator \cite{reg1}). For an algebraic number
field $K$ and prime $\ell$, let $\varepsilon_1,\ldots,\varepsilon_r$ be a system
of fundamental units of $K$ and $\varepsilon_0=1+\ell$ ($\varepsilon_0=5$ if
$\ell=2$). Put $R_\ell(K):=\det({\rm Sp}_{K/\mathbb
Q}(\log_\ell\varepsilon_i\cdot\log_\ell\varepsilon_j)),\quad 0\leqslant
i,j\leqslant r$, where by ${\rm Sp}_{K/\mathbb Q}$ we denote the mapping
$K\otimes_\mathbb Q \mathbb Q_\ell\to\mathbb Q_\ell$ induced by the trace
map ${\rm Sp}_{K/\mathbb Q}\colon K\to\mathbb Q$. Then $R_\ell(K)\neq 0$ for
any $K$ and $\ell$.
\label{conj4}
\end{conj}
A relation of the following conjecture to Conjecture~\ref{conj4} is the same as
that of the feeble Leopoldt conjecture and Leopoldt conjecture. But in
distinction with the feeble Leopoldt conjecture, the conjecture, which we
formulate below, is stated only in some particular cases. To formulate it, we
need a notion of the relative $\ell$-adic regulator. Let $K/k$ be an extension
of algebraic number fields and $U(K/k)$ the group of relative units of $K$ over
$k$, that is,
\[
 U(K/k)=\{\, x\in U(K)\mid N_{K/k}(x)\in\mu(k)\,\},
\]
where $\mu(k)$ is the group of all roots of unity in $k$, and $u_1,\ldots,u_t$
some system of fundamental units of $U(K/k)$. It means that $u_1,\ldots, u_t$
are independent and generate with the group $\mu(K)$ all the group $U(K/k)$.
Then the relative $\ell$-adic regulator $R_\ell(K/k)$ is defined by the formula
\[
 R_\ell(K/k)=\det\bigl({\rm Sp}_{K/\mathbb Q}(\log_\ell u_i\log_\ell
u_j)\bigr),\quad 1\leqslant i,j\leqslant t.
\]
\begin{conj}
 (Feeble conjecture on the $\ell$-adic regulator \cite{reg1}). For an algebraic
number field $K$ and its cyclotomic $\mathbb Z_\ell$-extension $K_\infty$, there
is an index $n_0$, which depends only on $K$ and $\ell$, such that
$R_\ell(K_m/K_n)\neq 0$ for any $m>n\geqslant n_0$.
\end{conj}
The regulator $\mathfrak R_\ell(K)$, which we determine and learn in the present
paper, characterize the group of universal global norms $U_{S,1}(K)$ defined
in~\eqref{1.0}. So we need a characterization of this group as Galois module. In
the next theorem we put no restrictions on decomposition type of prime $\ell$ in
 $K$.
\begin{theorem}
 Let $K$ be a finite Galois extension of $\mathbb Q$ with a Galois group $G$ and
$R=\mathbb Q_\ell[G]$. Suppose that the Leopoldt conjecture holds for
$(K,\ell)$. Then there is an isomorphism of $R$-modules
\[
 \varphi\colon U_{S,1}(K)\otimes_{\mathbb Z_\ell}\mathbb
Q_\ell\cong(U(K)\otimes_\mathbb Z \mathbb Q_\ell)\oplus\mathbb Q_\ell
\]
(here we write the group operation in $U_{S,1}(K)$ and $U(K)$ additively).
Without assumption that $K$ is Galois over $\mathbb Q$, the module $U_{S,1}(K)$
has $\mathbb Z_\ell$-rank $r_1+r_2$, where $r_1$ and $r_2$ is the number of real
and complex places of $K$ respectively.
\end{theorem}
{\bf Proof.} We shall state only the first part of the theorem. To state the
second one, man has to repeat the proof of the first part substituting
everywhere isomorphisms of $R$-modules by isomorphisms of corresponding $\mathbb
Q_\ell$-spaces.

 It follows from the representation theory of finite groups, that any finitely
generated $R$-module may be uniquely presented as a direct sum of irreducible
$R$-modules. In particular, for $R$-modules $A,B,C$ the condition $A\oplus
B\cong A\oplus C$ yields $B\cong C$.

Let $M$ be the maximal Abelian $\ell$-extension of $K_\infty$ unramified outside
$\ell$ and $X=G(M/K_\infty)$. Let $M_n$ be the maximal Abelian $\ell$-extension
of $K_n$ unramified outside $\ell$.
Then for any $n$ there are inclusions $K\subseteq K_n\subset K_\infty\subset
M_n\subseteq M$, and the Galois group $G(M_n/K_n)$ enters the exact sequence
\[
 1\longrightarrow X_n\longrightarrow
G(M_n/K_n)\longrightarrow\Gamma_n\longrightarrow 1,
\]
where $X_n=G(M_n/K_\infty)=X/(\gamma_n-1)X,\quad \gamma_n=\gamma^{\ell^n}$ and
$\Gamma_n=G(K_\infty/K_n)=\langle\gamma_n\rangle$. In particular, for $n=0$ one
gets an exact sequence of $\mathbb Z_\ell[G]$-modules
\begin{equation}
 1\longrightarrow X_0\longrightarrow
G(M_0/K)\longrightarrow\Gamma\longrightarrow 1.
\label{2.2}
\end{equation}
To prove the theorem, we shall determine the $R$-module
$Y(K):=G(M_0/K)\otimes_{\mathbb Z_\ell}\mathbb Q_\ell$ in two different ways.

Let $V(K)$ be a subgroup of $G(M_0/K)$ generated by the inertia subgroups for all
$v$ over $\ell$. Then the $G$-module $V(K)$ is of finite index in $G(M_0/K)$,
and, according to global class field theory, $V(K)$ contains in the exact
sequence of $G$-modules
\begin{equation}
 1\longrightarrow
U(K)[\ell]\stackrel{\alpha}{\longrightarrow}\prod_{v|\ell}
U(K_v)_\ell\longrightarrow V(K)\longrightarrow 1,
\label{2.3}
\end{equation}
where $U(K_v)$ is the group of units of the local field $K_v, \quad U(K_v)_\ell$
is the $\ell$-component of $U(K_v)$ (thus it coincides with the group of
principal units $U^{(1)}(K_v)$), and the product is taken over all places $v$ of
$K$ over $\ell$. Note that $\alpha$ is an injection, since we assume that the
Leopoldt conjecture holds for $(K,\ell)$.

Since the $\ell$-adic logarithm maps $U^{(1)}(K_v)$ on some full lattice in
$K_v$  and has finite kernel, it induces an isomorphism of $R$-modules
\begin{equation}
 \prod_{v|\ell}\bigl(U(K_v)_\ell\otimes_{\mathbb Z_\ell}\mathbb
Q_\ell\bigr)\cong\prod_{v|\ell}K_v\cong K\otimes_\mathbb Q\mathbb Q_\ell\cong R,
\label{2.4}
\end{equation}
where the last isomorphism follows from the theorem on normal base. Thus the
exact sequence~\eqref{2.3} yields an isomorphism
\begin{equation}
 \bigl(U(K)[\ell]\otimes_{\mathbb Z_\ell}\mathbb Q_\ell\bigr)\oplus Y(K)\cong R,
\label{2.5}
\end{equation}
where $Y(K):=V(K)\otimes_{\mathbb Z_\ell}\mathbb Q_\ell$.
Now we shall characterize the $G$-module $Y(K)$ in another way. To do this, we
calculate the $G$-module $X_0=G(M_0/K_\infty)$ in terms of Iwasawa theory. Note
that $G(M_n/K_n)$ enters the exact sequence of $\Gamma/\Gamma_n$-modules
\begin{equation}
 1\longrightarrow W(K_n)\longrightarrow G(M_n/K_n)\longrightarrow
C(K_n)\longrightarrow 1,
\label{2.6}
\end{equation}
where $W(K_n)$ is a subgroup of $G(M_n/K_n)$ generated by the decomposition
subgroups of all places $v$ over $\ell$ and $C(K_n)$ is a Galois group of the
maximal Abelian unramified $\ell$-extension of $K_n$, in which all places
$v|\ell$ splits completely. Note that by global class field theory $C(K_n)$
is canonically isomorphic to ${\rm Cl}'_\ell(K_n)$, where ${\rm Cl}'_\ell(K_n)$ is a factor group of
the $\ell$-component ${\rm Cl}_\ell(K_n)$ of the class group ${\rm Cl}(K_n)$ of
$K_n$ by a subgroup generated by all prime divisors of $\ell$.

For any pair of indices $m>n$ man can include the sequences \eqref{2.6} for
$K_n$ and $K_m$ into the commutative diagram
\begin{equation}
 \begin{CD}
  1@>>> W(K_m)@>>> G(M_m/K_m)@>>>C(K_m)@>>>1\\
@. @VVf_{m,n}V  @VVV @VVV @.\\
1 @>>> W(K_n)@>>> G(M_n/K_n)@>>>C(K_n)@>>>1
 \end{CD}
\label{2.6a}
\end{equation}
The vertical arrows of this diagram are induced by the restriction of
automorphisms.

Passing to projective limit in~\eqref{2.6} with respect to the vertical arrows,
we obtain an exact sequence of $\Gamma$-modules
\begin{equation}
 1\longrightarrow W(K_\infty)\longrightarrow X\longrightarrow
T_\ell(K_\infty)\longrightarrow 1,
\label{2.7}
\end{equation}
where $W(K_\infty)$ is a subgroup of $X$ generated by decomposition subgroups of
all $v$ over $\ell$ and $T_\ell(K_\infty)$ a Galois group of the maximal
Abelian unramified $\ell$-extension of $K_\infty$, in which all places of
$K_\infty$ over $\ell$ split completely.

By Proposition~\ref{prop2.2} one has $X^\Gamma=1$, hence \eqref{2.7} yields an exact
homological sequence
\begin{equation}
 1=X^\Gamma\longrightarrow T_\ell(K_\infty)^\Gamma\longrightarrow
W(K_\infty)_0\longrightarrow X_0\longrightarrow
T_\ell(K_\infty)_0\longrightarrow 1,
\label{2.8}
\end{equation}
where for a $\Gamma$-module $A$ we denote by $A_0$ the group $A/(\gamma-1)A$.

By Iwasawa theory $T_\ell(K_\infty)\otimes_{\mathbb Z_\ell}\mathbb Q_\ell$ is a
finite-dimensional $\mathbb Q_\ell$-space.
\begin{lemma}
 Let $K$ be a finite Galois extension of $\mathbb Q$ with Galois group $G$ and
$R=\mathbb Q_\ell[G]$. Let $A$ be a finite-dimensional $\mathbb Q_\ell$-space
with continuous action of the Galois group $G(K_\infty/\mathbb Q)$. Then the
$R$-modules $A^\Gamma$ and $A_0=A/(\gamma-1)A$ are isomorphic.
\end{lemma}
{\bf Proof.} We use induction by dimension of $A$. If $\dim A=1$, then the
statement of the lemma is obvious. Let $\dim A=d$ and assume that the lemma
holds for all dimensions less than $d$. If $A^\Gamma=0$, then multiplication by
$(\gamma-1)$ is an injection and thus an isomorphism. So $A_0=0$ and the lemma
holds. Suppose that $A^\Gamma\neq 0$. Put $B=A^\Gamma\cap (\gamma-1)A$. If $B=0$
then, taking into account that $\dim A^\Gamma=\dim A_0$, we see that the natural
surjection $A\to A_0$ induces an isomorphism of $R$-modules $A^\Gamma\to A_0$.
If $B\neq 0$ then we put $B_1=\{\, x\in A/B\mid (\gamma-1)\bar x\in
B\,\}=(A/B)^\Gamma$, where $\bar x$ is some pre-image of $x$ in $A$. Then
multiplication by $\gamma-1$ defines an $R$-surjection $B_1\to B$, whose kernel
coincides with $A^\Gamma/B$. Therefore, $B_1\cong A^\Gamma$ as $R$-modules, that
is, the $R$-modules $A^\Gamma$ and $(A/B)^\Gamma$ are isomorphic. Since $\dim
B\neq 0$, we have $(A/B)^\Gamma\cong (A/B)_0$ by assumption of induction but
$(A/B)_0\cong A_0$. This proves the lemma.

Taking tensor multiplication of~\eqref{2.8} by $\mathbb Q_\ell$ and applying
Lemma~2.1 to the $R$-module $T_\ell(K_\infty)\otimes_{\mathbb Z_\ell}\mathbb
Q_\ell$, we obtain
\begin{equation}
 W(K_\infty)_0\otimes_{\mathbb Z_\ell}\mathbb Q_\ell\cong X_0\otimes_{\mathbb
Z_\ell}\mathbb Q_\ell.
\label{2.9}
\end{equation}
By global class field theory for any $n\geqslant 0$ the group $W(K_n)$ enters
the exact sequence of $G(K_n/\mathbb Q)$-modules
\begin{equation}
1\lr E(K_n)\lr U_S(K_n)[\ell]\lr \prod_{v|\ell}(K_{n,v}^\times[\ell])\lr
W(K_n)\lr 1,
\label{2.10}
\end{equation}
where $[\ell]$ means pro-$\ell$-completion. If $m>n$ then man can include the
sequences~\eqref{2.10} for $K_m$ and $K_n$ in the following commutative diagram
$$
\begin{CD}
 1@>>>
E_\ell(K_m)@>>>U_S(K_m)[\ell]@>>>\prod_{v|\ell}(K_{m,v}^\times[\ell]
)@>>>W(K_m)@>>>1\\
@. @VV{N_{m,n}}V @VV{N_{m,n}}V @VV{N_{m,n}}V @VVf_{m,n}V @.\\
1@>>> E_\ell(K_n)@>>>
U_S(K_n)[\ell]@>>>\prod_{v|\ell}(K_{n,v}^\times[\ell])@>>>W(K_n)@>>> 1
\end{CD}
$$ 
where the maps $N_{m,n}$ are induced by the norm maps from $K_m$ into $K_n$ and
$f_{m,n}$ is the map of the diagram~\eqref{2.6a}.

Passing in~\eqref{2.10} to projective limit with respect to the vertical maps of
preceding diagram and taking into account that $\varprojlim E_\ell(K_n)=1$ by
Proposition~\ref{prop2.1}, we get an exact sequence of $\Gamma$-modules
\begin{equation}
 1\lr U_S(K_\infty)\lr\varprojlim\prod_{v|\ell}(K_{n,v}^\times[\ell])\lr
W(K_\infty)\lr 1,
\label{2.11}
\end{equation}
where $U_S(K_\infty)=\varprojlim U_S(K_n)[\ell]$ and the limit is taken with
respect to the norm maps. Denote the group
$\varprojlim\prod_{v|\ell}(K_{n,v}^\times[\ell])$ by $\tilde H(K_\infty)$.
Since $W(K_\infty)\subseteq X$, it follows from Proposition~\ref{prop2.2} that
$W(K_\infty)^\Gamma=0$, therefore, \eqref{2.11} induces an exact homological
sequence
\begin{equation}
 1\lr U_S(K_\infty)_0\lr \tilde H(K_\infty)_0\lr W(K_\infty)_0\lr 1.
\label{2.11a}
\end{equation}
Suppose that for any field $K_n$, where $n\geqslant 0$ one has a Galois $\mathbb
Z_\ell$-module $A(K_n)$. Moreover, suppose that for any pair of indices $m>n$ one has an
injection $i(n,m)\colon A(K_n)\to A(K_m)$ such that
$i(n,m)(A(K_n))=A(K_m)^{\Gamma_n}$, and for any triple of $m>n>k$ one has $i(k,m)=i(n,m)\circ
i(k,n)$. Then for any $m>n$ the norm map $N_{\Gamma_n/\Gamma_m}\colon A(K_m)\to
A(K_m)^{\Gamma_n}$ with respect to the $\Gamma_n/\Gamma_m$ induces a map
$N_{m,n}\colon A(K_m)\to A(K_n),\quad N_{m,n}=i(n,m)^{-1}\circ
N_{\Gamma_n/\Gamma_m}$. Put $A(K_\infty)=\varprojlim A(K_n)$, where the limit is
taken with respect to the maps $N_{m,n}$. Then the collection of the maps
$N_{m,0}\colon A(K_m)\to A(K)$ defines a map $N\colon A(K_\infty)\to A(K)$. On
the other hand, any map $N_{m,n}$ permits passing through the natural projection
$\pi_m\colon A(K_m)\to A(K_m)/(\gamma-1)A(K_m)$. Thus, the map $N$ permits
passing through the natural projection $\pi\colon A(K_\infty)\to
A(K_\infty)_0=A(K_\infty)/(\gamma-1)A(K_\infty)$, that is, there exists a map
$j\colon A(K_\infty)_0\to A(K)$ such that $N=j\circ \pi$. We shall deal with two
such collections of Galois modules: $A_1(K_n)=U_S(K_n)[\ell]$ and
$A_2(K_n)=\prod_{v|\ell}(K_{n,v}^\times[\ell])$. Denoting the maps $N,\pi$ and
$j$ for these collections by $N_1,\pi_1,j_1$ and $N_2,\pi_2,j_2$ respectively,
we see that man can include the exact sequences~\eqref{2.11} and~\eqref{2.11a}
into a commutative diagram with exact lines, where $j_1\circ\pi_1=N_1$ and
$j_2\circ\pi_2=N_2$,
\begin{equation}
 \begin{CD}
  1@>>> U_S(K_\infty) @>>>\tilde H(K_\infty) @>>> W(K_\infty) @>>> 1\\
@. @VV\pi_1V @VV\pi_2V @VV\pi_3V @.\\
1@>>> U_S(K_\infty)_0 @>\alpha>> \tilde H(K_\infty)_0 @>>> W(K_\infty)_0 @>>>
1\\
@. @VVj_1V @VVj_2V @. @. \\
{ } @. U_S(K)[\ell] @>\beta>> \prod_{v|\ell}(K_v^\times[\ell]) @. { } @. { } 
 \end{CD}
\label{2.12}
\end{equation}
\begin{lemma}
 The homomorphisms $j_1$ and $j_2$ are injections.
\end{lemma}
{\bf Proof.} At first, we prove that $j_2$ is an injection. To do this, it is
enough to check that $\hat
H^{-1}(\Gamma/\Gamma_n,\prod_{v|\ell}(K_{n,v}^\times[\ell]))=0$ for any $n$,
where $\hat H^{-1}$ is the Tate group of cohomologies. Since $K$ is normal over
$\mathbb Q$, there is an index $r\geqslant 0$ such that all places of $K$ over
$\ell$ split completely in $K_r/K$ and any place $v|\ell$ of $K_r$ has the
unique extension to $K_\infty$. If $n\leqslant r$ then the
$\Gamma/\Gamma_n$-module $D_n:=\prod_{v|\ell}(K_{n,v}^\times[\ell])$ is induced,
hence cohomologically trivial. If $n>r$ then $D_n$ is relatively induced.
Therefore by the Shapiro lemma one has $\hat H^{-1}(\Gamma/\Gamma_n,D_n)=\hat
H^{-1}(\Gamma_r/\Gamma_n,\prod_{v\in S_0}(K_{n,v}^\times[\ell]))=\prod_{v\in
S_0}\hat H^{-1}(\Gamma_r/\Gamma_n,K_{n,v}^\times[\ell])$, where $S_0$ is a set
of places of $K_n$, which for any place $v_0|\ell$ of $K$ contains exactly one
its extension to $K_n$. The last group is zero by virtue of Gilbert's Theorem
90, hence $\hat H^{-1}(\Gamma/\Gamma_n,D_n)=0$. It means that $\pi_2$ and $N_2$
have the same kernel, hence $j_2$ is an injection. Thus, $j_2$ settles an
isomorphism between $\tilde H(K_\infty)_0$ and the subgroup of (local) universal
norms in $\prod_{v|\ell}(K_v^\times[\ell])$.

Since $j_2\circ\alpha=\beta\circ j_1$ and $\alpha, \beta, j_2$ are injections,
we obtain that $j_1$ is also an injection. This proves the lemma.

Thus, $j_1$ settles an isomorphism between $U_S(K_\infty)_0$ and the group of
global universal norms ${\rm Im}N_1=U_{S,1}(K)$.

For a finite extension $L$ of $\mathbb Q_\ell$, let $F(L)$ be the subgroup of
all universal norms in $L^\times[\ell]$ from $L_\infty$. In other words,
$F(L)=\cap_{n=1}^\infty N_{L_n/L}(L_n^\times[\ell])$. Then, identifying in the
second line of~\eqref{2.12} the modules $U_S(K_\infty)_0$ and $\tilde
H(K_\infty)_0$ with $U_{S,1}(K)$ and $\prod_{v|\ell}F(K_v)$ (the identification
is given via the maps $j_1$ and $j_2$ respectively), we get an exact sequence
\begin{equation}
 1\lr U_{S,1}(K)\lr \prod_{v|\ell}F(K_v)\lr W(K_\infty)_0\lr 1.
\label{2.13}
\end{equation}
\begin{lemma}
 Let $K$ be a Galois extension of $\mathbb Q$ with Galois group $G$ and
$R=\mathbb Q_\ell[G]$. Then there is an isomorphism of $R$-modules
\[
 \bigl(\prod_{v|\ell}F(K_v)\bigr)\otimes_{\mathbb Z_\ell}\mathbb Q_\ell\cong R.
\]
\end{lemma}
{\bf Proof.} Let $L$ be a finite Galois extension of $\mathbb Q_\ell$ with
Galois group $H$ and $\theta_L\colon L^\times[\ell]\hookrightarrow
G_{\ell,L}^{ab}$ the reciprocity map, where $G_{\ell,L}^{ab}$ is a Galois group
of the maximal Abelian $\ell$-extension of $L$. Let $\pi\colon
G_{\ell,L}^{ab}\to\Gamma$ be the natural projection, where $\Gamma$ is a Galois
group of the cyclotomic $\mathbb Z_\ell$-extension of $L$. Then $F(L)$ coincides
with the kernel of the map $\pi\circ\theta_L$. If $L=\mathbb Q_\ell$ then one
can check immediately that $F(\mathbb Q_\ell)$ coincides with the kernel of the
logarithmic map $\log\colon\mathbb Q_\ell^\times[\ell]\to \mathbb Q_\ell$ (note
that the map $\log$, that was previously defined on $\mathbb Q_\ell^\times$, can
be uniquely extended on $\mathbb Q_\ell^\times[\ell]$ by linearity). 

If $L$ is a finite extension of $\mathbb Q_\ell$ then the following diagram is
commutative
$$
\begin{CD}
 L^\times[\ell]@>\theta_L>> G_{\ell,L}^{ab}\\
@VVN_{L/\mathbb Q_\ell}V @VVresV \\
\mathbb Q_\ell^\times[\ell]@>\theta_{\mathbb Q_\ell}>> G_{\ell,\mathbb
Q_\ell}^{ab}
\end{CD}
$$
where $N_{L/\mathbb Q_\ell}$ is the norm map and $res$ the restriction of
automorphisms. Thus, $F(L)$ coincides with pre-image of $F(\mathbb Q_\ell)$ with
respect to the map $N_{L/\mathbb Q_\ell}$. In other words,
\begin{equation}
 F(L)=\{\, x\in L^\times[\ell]\mid {\rm Sp}_{L/\mathbb Q_\ell}(\log x)=0\,\}.
\label{2.14}
\end{equation}
Put $P(L)=F(L)\cap U^{(1)}(L)$, where $U^{(1)}(L)$ is the group of principal
units of $L$. Then, for any $x\in P(L)$ man has $N_{L/\mathbb Q_\ell}(x)=1$ if
$\ell\neq2$ and $N_{L/\mathbb Q_\ell}(x)=\pm1$ if $\ell=2$. Hence the group
$U^{(1)}(L)$ enters the exact sequence of $H$-module
\[
 1\lr P(L)\lr U^{(1)}(L)\stackrel{\alpha}{\lr}\mathbb Z_\ell\lr 0,
\]
where $\alpha=\log\circ N_{L/\mathbb Q_\ell}$. Thus, using additive notation for
the group operation in $P(L)$ and $U^{(1)}(L)$, we get
\begin{equation}
 \bigl(P(L)\otimes_{\mathbb Z_\ell}\mathbb Q_\ell\bigr)\oplus \mathbb
Q_\ell\cong U^{(1)}(L)\otimes_{\mathbb Z_\ell}\mathbb Q_\ell\cong L\cong \mathbb
Q_\ell[H].
\label{2.15}
\end{equation}
On the other hand, there is an exact sequence of $H$-modules
\begin{equation}
 1\lr P(L)\lr F(L)\stackrel{\nu}{\lr} A\lr 1,
\label{2.16}
\end{equation}
where for $x\in F(L)$ by $\nu(x)\in \mathbb Z_\ell$ we denote the $\ell$-adic
exponent of $x$. The map $\nu$ is non-zero since $\ell\in F(L)$ and
$\nu(\ell)\neq0$, therefore $A$ is a subgroup of finite index in $\mathbb
Z_\ell$, hence $A\cong \mathbb Z_\ell$.
Comparing \eqref{2.15} and \eqref{2.16}, we obtain $F(L)\otimes_{\mathbb
Z_\ell}\mathbb Q_\ell\cong U^{(1)}(L)\otimes_{\mathbb Z_\ell}\mathbb
Q_\ell\cong\mathbb Q_\ell[H]$.

Now we return to $K$. As it was proved before, for any its completion $K_v$,
where $v|\ell$, one has $F(K_v)\otimes_{\mathbb Z_\ell}\mathbb
Q_\ell\cong\mathbb Q_\ell[G_v]$, where $G_v$ is a decomposition subgroup of
 $v$ in $G$. Since $\bigl(\prod_{v|\ell}F(K_v)\bigr)\otimes_{\mathbb
Z_\ell}\mathbb Q_\ell\cong \prod_{v|\ell}\bigl(F(K_v)\otimes_{\mathbb
Z_\ell}\mathbb Q_\ell\bigr)$ and $\prod_{v|\ell}F(K_v)\cong {\rm
Ind}_G^{G_v}F(K_v)$, we obtain
\[
 \bigl(\prod_{v|\ell}F(K_v)\bigr)\otimes_{\mathbb Z_\ell}\mathbb Q_\ell\cong
{\rm Ind}_G^{G_v}\bigl(F(K_v)\otimes_{\mathbb Z_\ell}\mathbb
Q_\ell\bigr)\cong{\rm Ind}_G^{G_v}\mathbb Q_\ell[G_v]\cong R.
\]
This proves the lemma.

Now we can easily conclude the proof of the theorem. It follows from \eqref{2.2}
and \eqref{2.5} that
\[
 (U(K)[\ell]\otimes_{\mathbb Z_\ell}\mathbb Q_\ell)\oplus(X_0\otimes_{\mathbb
Z_\ell}\mathbb Q_\ell)\oplus\mathbb Q_\ell\cong R.
\]
On the other hand, by \eqref{2.9}, \eqref{2.13} and Lemma~2.3 we get
\[
 (U_{S,1}(K)\otimes_{\mathbb Z_\ell}\mathbb Q_\ell)\oplus(X_0\otimes_{\mathbb
Z_\ell}\mathbb Q_\ell)\cong R.
\]
This proves the theorem.
\begin{con}
 Let $k=\mathbb Q(\sqrt{-d})$ be an imaginary quadratic field and $K$ a Galois
extension of $\mathbb Q$ with a Galois group $G$, such that $K$ contains $k$.
Put $H=G(K/k)$ and let $\tau$ be an automorphism of complex conjugation in $G$.
In particular, it means that $\tau^2=1$ и $\tau(\sqrt{-d})\neq\sqrt{-d}$. Then
in the group $U_{S,1}(K)$ there is an element $\omega$ such that
$\tau(\omega)=\omega$ and $\{h(\omega)\}_{h\in H}$ generate a submodule of
finite index in $U_{S,1}(K)$.
\end{con}
Indeed, in the group of units $U(K)$ there is the so cold an Artin unit
$\varepsilon$ such that $\tau(\varepsilon)=\varepsilon$ and
$\{h(\varepsilon)\}_{h\in H}$ generate in $U(K)$ a subgroup of finite index.
Then $\varepsilon_1=\varepsilon\oplus 1$ generates an $R$-module
$\bigl(U(K)\otimes_\mathbb Z\mathbb Q_\ell\bigr)\oplus\mathbb Q_\ell$. So the
element $\omega_1=\varphi^{-1}(\varepsilon_1)$, where $\varphi$ is an
isomorphism of Theorem 1, has the property $\tau(\omega_1)=\omega_1$. Then the
elements $\{h(\omega_1)\}_{h\in H}$ form the basis of $\mathbb Q_\ell$-space
$U_{S,1}(K)\otimes_{\mathbb Z_\ell}\mathbb Q_\ell$. If $s$ is sufficiently large
then $\ell^s\omega_1\in U_{S,1}(K)$, so one can put $\omega=\ell^s\omega_1$.
\section{Conditional proof of the main conjecture}
\setcounter{equation}{0}
Firstly, we give more detailed definition of the regulator $\mathfrak R_\ell(K)$
for an algebraic number field $K$ and prime $\ell$ that splits completely
in $K$.

Let $\mathfrak l_1,\ldots,\mathfrak l_m$ run all prime divisors over $\ell$ in
$K$. Then for any $S$-unit $u\in U_S(K)$, where $S=\{\mathfrak
l_1,\ldots,\mathfrak l_m\}$, there is its principal divisor $(u)$, which we
shall write in additive form as $(u)=\sum_{i=1}^ma_i\mathfrak l_i$, where
$a_i\in\mathbb Z$. We define the group $D(K)$ as a free $\mathbb Z_\ell$-module
with generators $\mathfrak l_1,\ldots,\mathfrak l_m$ and shall consider $(u)$ as
an element of $D(K)$. Then to any $u$ from the pro-$\ell$-com\-ple\-tion $U_S(K)[\ell]$
of the group $U_S(K)$ corresponds its divisor $\di(u)\in D(K)$. If
$u=\lim_{n\to\infty}u_n$, where $u_n\in U_S(K)$, then
$\di(u)=\lim_{n\to\infty}(u_n)$. Obviously, $\di(u)$ depends only on the image
of $u$ in $\bar U_S(K)[\ell]:=U_S(K)[\ell]/\mu_\ell(K)$, where $\mu_\ell(K)$ is
the group of all roots of unity of $\ell$-power degree in $K$.

We consider on the $\mathbb Z_\ell$-module $D(K)$ a symmetric non-degenerate
bilinear form defined for  $x,y\in D(K),\quad x=\sum_{i=1}^ma_i\mathfrak
l_i,\quad y=\sum_{i=1}^mb_i\mathfrak l_i$ by
\begin{equation}
 \langle x,y\rangle=\sum_{i=1}^m a_ib_i\in\mathbb Z_\ell.
\label{3.1}
\end{equation}
If we wish to stress dependence of the product~\eqref{3.1} upon $K$ we shall
write $\langle x,y\rangle_K$.

If $\sigma$ is an automorphism of $K$ then $\sigma$ acts on $S$ and hence on the
group $D(K)$. Obviously, one has
\begin{equation}
 \langle x,y\rangle=\langle\sigma x,\sigma y\rangle.
\label{3.2}
\end{equation}
Let $K_\infty$ be the cyclotomic $\mathbb Z_\ell$-extension of $K$ and
$U_{S,1}(K)$ the group of global universal norms defined in the introduction.
Let $e_1,\ldots,e_t$ be a $\mathbb Z_\ell$-basis of $U_{S,1}(K)$. Then we put
\[
 \mathfrak R_\ell(K)=\det\bigl(\langle\di(e_i),\di(e_j)\rangle\bigr),\qquad
1\leqslant i,j\leqslant t.
\]
If $e'_1,\ldots,e'_t$ is another basis of $U_{S,1}(K)$ and $C$ is a transition
matrix then the regulators defined by this two basis differs by
factor $(\det C)^2\in \mathbb Z_\ell^\times$. In particular, the $\ell$-adic
exponent $\nu_\ell\bigl(\mathfrak R_\ell(K)\bigr)$ does not depend on the choice
of the basis in $U_{S,1}(K)$.

Note some simple properties of pairing $\langle\phantom x,\phantom y\rangle$ and
the regulator $\mathfrak R_\ell(K)$.
\begin{prop}
Suppose that $\ell$ splits completely in $K$, and Leopoldt conjecture holds
for $(K,\ell)$. 
Then the map $\di\colon U_{S,1}(K)\to D(K)$ is an injection.
 \label{prop3.1}
\end{prop}
{\bf Proof.} The kernel of the map $\di\colon \bar U_S(K)[\ell]\to D(K)$
coincides with the pro-$\ell$-completion of $\bar U(K)[\ell]$, where $\bar
U_S(K)=U_S(K)/\mu(K)$, $\bar U(K)=U(K)/\mu(K)$ and $\mu(K)$ is the group of all
roots of unity in $K$. Thus, it is enough to check that $U_{S,1}(K)\cap(\bar
U(K)[\ell])=1$. If $\xi$ is an element of this intersection, then, as was shown
in the proof of Theorem~1, for any injection $\sigma\colon K\hookrightarrow
\mathbb Q_\ell$ man has $\log\sigma(\xi)=0$. It means that $\xi^d$, where $d$ is
an order of $\mu(K)$, belongs to the Leopoldt kernel $E_\ell(K)$, which vanishes
by assumption. This proves the proposition.

Let $L/K$ be an extension of algebraic number fields and prime $\ell$ splits
completely in $L$, hence also in $K$. Then there are natural maps $i_{L/K}\colon
D(K)\to D(L)$ and $N_{L/K}\colon D(L)\to D(K)$. The first one is induced by the
injection $K\hookrightarrow L$ and the second one by the norm map from $L^\times$ to
$K^\times$. Then for any $x,y\in D(K)$ и $z\in D(L)$ man has
\begin{equation}
 \langle x,y\rangle_K=[L:K]\langle i_{L/K}(x),i_{L/K}(y)\rangle_L,\quad \langle
x,N_{L/K}(z)\rangle_K=\langle i_{L/K}(x),z\rangle_L.
\label{3.3}
\end{equation}
Indeed, it is enough to check the equalities~\eqref{3.3} in the case, when
$x,y,z$ are prime divisors, but, obviously, in this case they hold true.
\begin{prop}
Let $L$ be a Galois extension of $\mathbb Q$ with Galois group $G$, prime $\ell$
splits completely in $L$ and $K\subseteq L$. Then the condition $\mathfrak
R_\ell(L)\neq0$ yields that $\mathfrak R_\ell(K)\neq0$.
\label{prop3.2}
\end{prop}
{Proof.} Put $H=G(L/K)$. The inclusion $K\hookrightarrow L$ induces the
inclusion $i\colon U_{S,1}(K)\hookrightarrow U_{S,1}(L)$, and the norm map
$N\colon L^\times\to K^\times$ induces the map $N\colon U_{S,1}(L)\to
U_{S,1}(K)$, whose kernel we denote by $B$. Since the group
$F:=i(U_{S,1}(K))\oplus B$ has finite index in $U_{S,1}(L)$ and $\mathfrak
R_\ell(L)\neq 0$ by assumption, we obtain that the pairing $\langle\phantom
x,\phantom x\rangle$ is non-degenerate on $F$. By~\eqref{3.3} the groups
$i(U_{S,1}(K))$ and $B$ are orthogolal with respect to this pairing, hence the
pairing $\langle\phantom x,\phantom X\rangle_L$ is non-degenerate on the group
$i(U_{S,1}(K))$. Then by~\eqref{3.3} the pairing $\langle\phantom x,\phantom
x\rangle_K$ is non-degenerate on $U_{S,1}(K)$, that is, $\mathfrak R_\ell(K)\neq
0$. This proves the proposition.

Though the definition of the regulator $\mathfrak R_\ell(K)$ does not use the
$\ell$-adic logarithms, we need them in order to give conditional proof of the
main conjecture (Theorem~2), so we need the following considerations.
\begin{conj}
 Let $K$ be a Galois extension of degree $n$ of $\mathbb Q$ with Galois group
$G$, and prime $\ell$ splits completely in $K$. Let $\mathfrak l$ be a prime
divisor of $K$ over $\ell$ and $\alpha$ such an element of $K^\times$ that
$(\alpha)=\mathfrak l^{\bf h}$ for some natural number ${\bf h}$ (for example,
man can take for $\bf h$ the class number of $K$). Then any $n-1$ of the elements
$\log_\ell\sigma_1(\alpha),\ldots,\log_\ell \sigma_n(\alpha)$, where $\sigma_i$
runs $G$, are lineary independent over $\mathbb Q_\ell$ in the space
$K\otimes_\mathbb Q\mathbb Q_\ell\cong \mathbb Q_\ell^n$.
\end{conj}
{\bf Remark.} Since $N_{K/\mathbb Q}(\alpha)=\ell^h$, one has a relation
\begin{equation}
 \sum_{i=1}^n\log_\ell\sigma_i(\alpha)=\log_\ell(N_{K/\mathbb Q}(\alpha))=0.
\label{3.3a}
\end{equation}
The next result means that Conjecture~5 is in some sense an analog of
Leopoldt conjecture. We will not use it in the present paper, but we give it for
completeness.
\begin{prop}
 Let $K$ be an Abelian algebraic number field of degree $n$ and prime $\ell$
splits completely in $K$. Then Conjecture~5 holds for $K$.
\label{prop3.3}
\end{prop}
{\bf Proof.} Put $G=G(K/\mathbb Q)$, and let $\alpha\in K^\times$ be as in
Conjecture~5. Put $\varepsilon_0=1+\ell$ if $\ell\neq2$ and $\varepsilon_0=5$ if
$\ell=2$. Put $\beta=\alpha\varepsilon_0$. Then to prove the proposition, it is
enough to check that an $n\times n$-matrix $A$, whose lines are the vectors
$\log_\ell\sigma_j\sigma_1^{-1}(\beta),\ldots,\log_\ell\sigma_j\sigma_n^{-1}
(\beta)$, is non-degenerate. A general element of this matrix is of the form
$a_{ij}=\log \sigma_j\sigma_i^{-1}(\beta)$. Therefore, the determinant $\det A$
of this matrix is Dedekind determinant, thus by \cite[Ch. 3, \S 6,
Theorem~6.1]{L} one has
\begin{equation}
 \det A=\prod_{\chi\in \hat G}\sum_{\sigma\in
G}\bar\chi(\sigma)\log\sigma(\beta),
\label{3.4}
\end{equation}
where $\hat G$ is the group of characters of $G$.

Obviously, the numbers $\sigma^{-1}_1(\beta),\ldots,\sigma^{-1}_n(\beta),\ell$
are multiplicatively independent. Hence, according to the $\ell$-adic analog of
Baker theorem \cite{B}, the numbers $\log
\sigma^{-1}_1(\beta),\ldots,\log\sigma^{-1}_n(\beta)\in \bar{\mathbb Q}_\ell$
are lineary independent over the field of all algebraic numbers $\bar{\mathbb
Q}$. Therefore, all the factors in the right hand part of \eqref{3.4} are
non-zero and $\det A\neq0$. This proves the proposition.
\begin{prop}
 Conjecture 5 follows from Conjecture 1.
\label{prop3.4}
\end{prop}
{\bf Proof.} Let $K/\mathbb Q$ be a Galois extension of degree $n$ with Galois
group $G$ and $\ell$ a prime number that splits completely in $K$. Let
$\varepsilon_0$ and $\beta$ have the same meaning as in the proof of
Proposition~\ref{3.3}. Reasoning as in the proof of that proposition, we see
that it is enough to state non-degenerateness of the $n\times n$-matrix with a
general element $a_{ij}=\log\sigma_j\sigma_i^{-1}(\beta),\quad 1\leqslant
i,j\leqslant n$. The determinant $\det A$ of this matrix is a polynomial of
degree $n$ in $n$ indeterminates
$X_1=\log\sigma_1(\beta),X_2=\log\sigma_2(\beta),\ldots,X_n=\log\sigma_n(\beta),
\,\, \sigma_i\in G$, where, without restriction of generality, one can assume
that $\sigma_1=1$. Then the matrix $A$ has the elements $X_1$ on the main
diagonal, and beyond it there are the elements $X_2,\ldots,X_n$. Thus, $\det
A=P(X_1,\ldots,X_n)=X_1^n+\ldots$, where $P$ is a polynomial of degree $n$ in
$X_1,\ldots,X_n$ and dots mean the sum of monomials that contain $X_1$ in degree
at most $n-1$. Hence $P$ is a non-zero polynomial. By Conjecture~1 the numbers
$\log\sigma_1(\beta),\ldots,\log\sigma_n(\beta)$ are algebraically independent,
thus $\det A\neq0$. This proves the proposition.

Together with the group of global universal norms $U_{S,1}(K)$ we consider also
the group of local universal  norms $U_{S,2}(K)$ defined by
\begin{equation}
 U_{S,2}(K):=\{\, x\in \bar U_S(K)[\ell]\mid x\in N_{n,v}(K_{n,v}^\times) \text{
for all } v\,\}.
\label{3.5}
\end{equation}
where $N_{n,v}$ is the norm map from completion $K_{n,v}$ of $K_n$ into $K_v$
and $v$ runs all the places of $K$. Note that $N_{n,v}$ is the norm map with
respect to the decomposition subgroup of $v$. 
It follows immediately from the definition that $U_{S,2}(K)\supseteq
U_{S,1}(K)$. 

Note that $K_n/K$ is unramified in any $v\nmid\ell$ and any $u\in U_S(K)$ is a
unit outside $S$. Therefore, the condition~\eqref{3.5} holds automatically for
all $v\nmid \ell$. If $v|\ell$ then by~\eqref{2.14} the condition~\eqref{3.5}
means that $\log(N_{K_v/\mathbb Q_\ell}(x))=0$. If $\ell$ splits completely
in $K$ then $K_v=\mathbb Q_\ell$ for any $v|\ell$. Thus~\eqref{3.5} is
equivalent to the condition $\log\sigma(x)=0$ for any inclusion of $K$ into
$\mathbb Q_\ell$, which, in turn, is equivalent to
\begin{equation}
 \log_\ell x=0.
\label{3.5a}
\end{equation}
\begin{prop}
 Let $K/\mathbb Q$ be a Galois extension of degree $n$ with Galois group $G$ and
prime $\ell$ splits completely in $K$. Assume that the pair $(K,\ell)$
satisfies Conjecture~2 (the Leopoldt conjecture) and Conjecture~5. Then the
group $U_{S,1}(K)$ has finite index in $U_{S,2}(K)$.
\label{prop3.5}
\end{prop}
{\bf Proof.} By Theorem~1 the $\mathbb Z_\ell$-rank of $U_{S,1}(K)$ equals
$r+1$, where $r$ is the rank of the group of units $U(K)$, so it is enough to
check that the same $\mathbb Z_\ell$-rank has $U_{S,2}(K)$.

Consider a $\mathbb Q_\ell$-space of dimension $n\quad K\otimes_\mathbb Q
\mathbb Q_\ell\cong \mathbb Q_\ell^n$. Let $V$ be a subspace of $\mathbb
Q_\ell^n$ generated by $\log_\ell\sigma_i(\alpha)$, where
$\sigma_1,\ldots,\sigma_n$ runs all inclusions of $K$ in $\mathbb Q_\ell$ and
$\alpha$ has the same meaning, as in Conjecture~5. Then any $x\in V$ may be
presented as a linear combination
\begin{equation}
 x=\sum_{i=1}^na_i\log_\ell\sigma_i(\alpha),
\label{3.6}
\end{equation}
where by~\eqref{3.3a} and Conjecture~5 the coefficients $a_i\in\mathbb Q_\ell$ are
defined uniquely up to a constant summand. In particular, in~\eqref{3.6} we may
assume that $a_n=0$ and $a_1,\ldots,a_{n-1}$ are uniquely defined by $x$.

Take for $x$ the logarithms of a system of fundamental units of $K$, that is,
$\log_\ell\varepsilon_1,\ldots,\log_\ell\varepsilon_r$. Thus, we obtain a system
of equalities
\[
 \sum_{i=1}^{n-1}c_{ij}\log_\ell\sigma_i(\alpha)=\log_\ell\varepsilon_j,\quad
1\leqslant j\leqslant r,
\]
where the coefficients $c_{ij}\in\mathbb Q_\ell$ are defined uniquely.

Let $\ell^s$ be such a power of $\ell$ that $d_{ij}:=\ell^sc_{ij}\in\mathbb
Z_\ell$ for all $i,j$. Put
\begin{equation}
 \eta_j=\prod_{i=1}^{n-1}\sigma_i(\alpha)^{d_{ij}},\quad 1\leqslant j\leqslant
r.
\label{3.6a}
\end{equation}
Then the elements $\eta_j\in U_S(K)[\ell]$ satisfy the condition
$\log_\ell\eta_j=0$, hence the images of the elements $\eta_j$ in
$U_S(K)[\ell]/\mu_\ell(K)$ belong to the group $U_{S,2}(K)$. Now suppose that
$y\in U_S(K)[\ell]$ and $y\bmod{\mu_\ell(K)}\in U_{S,2}(K)$. Then, putting
$y_1=y^{\ell^s}$ for sufficiently large $s$, we obtain a presentation
\begin{equation}
 y_1=\prod_{i=1}^n\sigma_i(\alpha)^{a_i}\times\prod_{j=1}^r\varepsilon_j^{b_j}
\label{3.7}
\end{equation}
for some $a_i,b_j\in\mathbb Z_\ell$. Since $\ell\in U_{S,2}(K)$, multiplying
$y_1$ by $\ell^{-a_n}$ and noting that $(\ell)=\mathfrak l_1\cdots\mathfrak
l_n$, we can assume that $a_n=0$ in~\eqref{3.7}.

Consider an element $z:=y_1\prod_{j=1}^r\eta_j^{b_j}$. Since $y_1\in
U_{S,2}(K)$, we have $\log_\ell y_1=0$, hence $\log_\ell z=0$. But~\eqref{3.6a}
and~\eqref{3.7} imply that $z=\prod_{i=1}^{n-1}\sigma_i(\alpha)^{h_i}$, где
$h_i=a_i+\sum_{j=1}^r b_j d_{ij}$. By Conjecture~5 the vectors
$\log_\ell\sigma_i(\alpha),\quad 1\leqslant i< n$, are lineary independent over
$\mathbb Q_\ell$. Therefore, $h_i=0$ for $i=1,\ldots,n-1$, that is, any $y\in
U_{S,2}(K)$ has a presentation of the form
$y=\ell^{b_0}\prod_{j=1}^r\eta_j^{b_j}\bmod{\mu_\ell(K)}$ for some
$b_0,b_1,\ldots,b_r\in\mathbb Q_\ell$. This proves that the $\mathbb
Z_\ell$-rank of $U_{S,2}(K)$ is at most $r+1$. Then Theorem~1 and the inclusion
$U_{S,1}(K)\subseteq U_{S,2}(K)$ yield that the index $(U_{S,2}(K):U_{S,1}(K))$
is finite. This proves the proposition.

{\bf Remark 1.} Let $K$ be any (maybe non-Galois over $\mathbb Q$) algebraic
number field and prime $\ell$ splits completely in $K$. Let $L$ be the Galois
closure of $K$. By Proposition~\ref{prop3.5} one has
$(U_{S,2}(L):U_{S,1}(L))<\infty$ if the conjectures~2 and~5 hold, but then these
conjectures yields that $(U_{S,2}(K):U_{S,1}(K))<\infty$.

{\bf Remark 2.} Let $T_\ell(K_\infty)$ be a $\Gamma$-module that enters the
proof of Theorem~1. In \cite{mt}, there was formulated a conjecture
(Conjecture~2) that $T_\ell(K_\infty)^\Gamma$ is finite for any algebraic number
field $K$ and any prime $\ell$. It was stated there (Proposition~7.5) that (in
notations of the present paper) $T_\ell(K_\infty)^\Gamma\cong
U_{S,2}(K)/U_{S,1}(K)$. Thus, if $\ell$ splits completely in $K$ then it
follows from Proposition~3.5 that the last conjecture is a consequence of
Conjecture~1.
\begin{prop}
 Assume that prime $\ell$ splits completely in $K$ and Leopoldt conjecture
holds for $(K,\ell)$. Then the group $U_{S,1}(K)$ has finite index in
$U_{S,2}(K)$ if and only if the image of the map $\log_\ell\colon
U_S(K)[\ell]\to\prod_{v|\ell}K_v$ is of $\mathbb Z_\ell$-rank $n-1$, where
$n=[K:\mathbb Q]$.
\label{prop3.6}
\end{prop}
{\bf Proof.} Let $\mathfrak l_1,\ldots,\mathfrak l_n$ be all prime divisors of
$\ell$ in $K$ and $\alpha_1,\ldots,\alpha_n\in U_S(K)$ such elements that
$(\alpha_i)=\mathfrak l^{s_i}_i$ for some natural $s_i$ и $i=1,2,\ldots,n$.
Let $N\subseteq\mathbb Z_\ell^n$ be a submodule of all vectors
$x=\{x_1,\ldots,x_n\}$ such that $\sum_{i=1}^n x_i\log_\ell\alpha_i=\log_\ell u$
for some $u\in U(K)[\ell]$. Then for any $x\in N$ the element
$u^{-1}\prod_{i=1}^n\alpha_i^{x_i}$ belongs to $U_{S,2}(K)$. Therefore, the rank
of $N$ equals to the rank of $U_{S,2}(K)$. If $U_{S,1}(K)$ has a finite index in
$U_{S,2}(K)$ then by Theorem~1 this rank equals $r_1+r_2$. It means that the
module $\mathbb Z_\ell^n/N$ is of the rank $r_2$.
So, a $\mathbb Z_\ell$-module $\log_\ell(U_S(K)[\ell])/\log_\ell(U(K)[\ell])$
has the rank $r_2$. Since a $\mathbb Z_\ell$-module $\log_\ell(U(K)[\ell])$ is
of the rank $r_1+r_2-1$, we have proved one half of the proposition.

Now we suppose that the $\mathbb Z_\ell$-rank of $\log_\ell(U_S(K)[\ell])$ is
$n-1$. It means that $\log_\ell (U_S(K)[\ell]) /\log_\ell(U(K)[\ell])$ has
$\mathbb Z_\ell$-rank $r_2$. But then the rank of $N$ is $r_1+r_2$, that is, the
rank of $U_{S,2}(K)$ coincides with the rank of $U_{S,1}(K)$. Hence the index 
$(U_{S,2}(K):U_{S,1}(K))$ is finite. This proves the proposition.

Now we can formulate and prove the main result of this section.
\begin{theorem}
 Let $K$ be an algebraic number field and prime $\ell$ splits completely in
$K$. Then Conjecture~1 yields that $\mathfrak R_\ell(K)\neq0$.
\end{theorem}
{\bf Proof.} If $\ell$ splits completely in $K$ then $\ell$ splits
completely in the normal closure $L$ of $K$. If $\ell$ splits in an
imaginary quadratic field $k=\mathbb Q(\sqrt{-d})$ then $\ell$ splits
completely in the field $F=L\cdot k$. By Proposition~\ref{prop3.2} the validity
of the theorem for $F$ yields its validity for $K$.

So, we can assume that $K$ is a Galois extension of $\mathbb Q$ with Galois
group $G$ and $K\supset k=\mathbb Q(\sqrt{-d})$. Let $H=G(K/k),\quad [K:k]=n$
and $\tau\in G$ be an automorphism of complex conjugation. In particular, this
means that $\tau^2=1$ и $\tau\notin H$. Let $\alpha\in U_S(K)$ be of the same
meaning, as in Conjecture~5. By Proposition~\ref{prop3.4} this conjecture holds
for $K$. Hence any  $2n-1$ of the $2n$ elements
$\{\log_\ell\sigma(\alpha)\}_{\sigma\in G}$ are lineary independent over
$\mathbb Q_\ell$ and form a basis of the $\mathbb Q_\ell$-space
\[
 V:=\{\,x\in K\otimes_\mathbb Q\mathbb Q_\ell\cong\mathbb Q_\ell^{2n}\mid {\rm
Sp}_{K/\mathbb Q}(x)=0\,\}.
\]
Let $\varepsilon\in U(K)$ be an Artin unit, that is,
$\tau(\varepsilon)=\varepsilon$ and the collection of $n$ elements
$\{\sigma(\varepsilon)\}_{\sigma\in H}$ generates the subgroup of finite index
in $U(K)$.  Since $\log_\ell(\varepsilon)\in V$, there are coefficients
$a_\sigma\in \mathbb Q_\ell,\quad\sigma\in G$, such that 
\[
 \sum_{\sigma\in G}a_\sigma\log_\ell\sigma(\alpha)=\log_\ell \varepsilon.
\]
Note that the coefficients $a_\sigma$ are defined uniquely up to an arbitrary
constant summand. That is, instead of $a_\sigma$ one can take the coefficients
$a'_\sigma=a_\sigma+c$, where $c\in\mathbb Q_\ell$ is independent of $\sigma$.

Put (in additive notation)
\[
 z=(1+\tau)\sum_{\sigma\in G}a_\sigma\sigma(\alpha)\in U_S(K)[\ell].
\]
Taking into account that $(1+\tau)\tau\sigma(\alpha)=(1+\tau)\sigma(\alpha)$, we
obtain
\begin{equation}
 z=(1+\tau)\sum_{h\in H}c_hh(\alpha),\qquad c_h=a_h+a_{\tau h}.
\label{3.7a}
\end{equation}
Here the coefficients $c_h$ are defined up to a constant summand and
\begin{equation}
 2\log_\ell\varepsilon=\log_\ell z=(1+\tau)\sum_{h\in H}c_h\log_\ell h(\alpha).
\label{3.8}
\end{equation}
Fix an element $h_0\in H,\quad h_0\neq 1$. Subtracting out of all $c_h$ the
coefficient $c_{h_0}$, we may assume, without restriction of generality, that
the coefficient $c_{h_0}$ in \eqref{3.8} is zero.
\begin{lemma}
 Suppose that $c_{h_0}=0$ in~\eqref{3.8} for some $h_0\neq1$. Then, assuming the
validity of Conjecture~1, man gets that $n-1$ coefficients $c_h$, $h\in
H$, $h\neq h_0$, are algebraically independent. 
\label{lemma3.1} 
\end{lemma}
{\bf Proof.} Consider a system of $3n-2$ elements $T=T_1\cup T_2$, where
$T_1=\{\sigma(\alpha)\}_{\sigma\in G,\sigma\neq 1},\quad
T_2=\{h(\varepsilon)\}_{h\in H,h\neq 1}$. Note that we consider $K$ as a
subfield of the field $\bar{\mathbb Q}_\ell$, and all logarithms that we
introduce below are the elements of this field. The system $T\cup\{\ell\}$ is
multiplicatively independent hence the logarithms of the elements of $T$, that
is, $3n-2$ elements $\log\sigma(\alpha)$, where $\sigma\in G$,  $\sigma\neq1$,
and $\log h(\varepsilon)$, where $h\in H$ and $h\neq1$, are algebraically
independent. In other words, the field generated by these logarithms, which we
denote by $\mathbb Q(\log T)$, is of transcendence degree $3n-2$ over $\mathbb
Q$.

Multiplying both parts of~\eqref{3.8} by any $\sigma\in H$ and taking into
account that the map $\log_\ell$ is a $G$-homomorphism, we get, using an
additive notation for multiplication,
\[
 2\log_\ell\sigma(\varepsilon)=\sum_{h\in H}c_h\log_\ell
\sigma(1+\tau)h(\alpha),\quad c_{h_0}=0.
\]
Equating the first coordinates in this equality of vectors, we obtain
\[
 2\log\sigma(\varepsilon)=\sum_{h\in H}c_h\log\sigma(1+\tau)h(\alpha)
\]
for all $\sigma\in H$. Thus, putting $T_3=\{c_\sigma\}_{\sigma\in H,\,\sigma\neq
h_0}$ (remind that $c_{h_0}=0)$, we obtain that all the elements of $\log T_2$
belong to the field $\mathbb Q(\log T_1,T_3)$. Since the transcendence degree of
$\mathbb Q(\log(T_1\cup T_2))$ is $3n-2$, we get that the transcendence degree
of $\mathbb Q(\log T_1,T_3)$ is at least $3n-2$. But $T_1\cup T_3$ is a set in
$3n-2$ elements. Therefore, all the elements of $T_3$ are algebraically
independent. This proves the lemma.

Continue the proof of the theorem. Put $z'=z\varepsilon^{-2}$. Then $\log_\ell
z'=0$, that is, $z'$ satisfy the condition~\eqref{3.5a}, hence also the
condition~\eqref{3.5}, that is, $z'\in U_{S,2}(K)$. Since we assume the
conjecture 1 to be true, we can assume that Leopoldt conjecture and Conjecture~5
hold for $K$.  Then by Proposition~\ref{prop3.5} we have $z_1:=(z')^{\ell^s}\in
U_{S,1}(K)$ for some natural $s$. We shall prove that the matrix $C=(\{h_i(z_1),
h_j(z_1)\}), \quad 1\leqslant i,j\leqslant n$
is non-degenerate, where $h_1,\ldots,h_n$ are all the elements of $H$ in some
order and $\{\phantom x,\phantom x\}$ means the product~\eqref{1.2}. If we shall
prove it then we get that the elements $h_1(z_1),\ldots,h_n(z_1)$ are
multiplicatively independent over $\mathbb Z_\ell$ (since their divisors are
lineary independent over $\mathbb Q_\ell$). Thus, by Theorem~1 the elements
$h_1(z_1),\ldots,h_n(z_1)$ generate a subgroup $Y$ of finite index in
$U_{S,1}(K)$. So, the condition $\det C\neq 0$ means that the pairing
$\{\phantom x,\phantom x\}$ is non-degenerate on $Y$. Then it is non-degenerate
on $U_{S,1}(K)$ as well. This will prove the theorem.

To simplify the notations, we shall instead of $C$ to have deal with the matrix
$C':=\ell^{-2s}\mathbf h^{-2}C$, where $\mathbf h$ is defined by condition
$(\alpha)=\mathfrak l^{\mathbf h}$. The general element of this matrix is
$\mathbf h^{-2}\{h_i(z'),h_j(z')\}$. Taking into account that
$\di(z')=\di(z)=\mathbf h(1+\tau)\sum_{\sigma\in H}c_\sigma\sigma(\mathfrak l)$,
we obtain
\begin{equation}
 \mathbf h^{-2}\{h_i(z'),h_j(z')\}=\langle h_i(1+\tau)\sum_{\sigma\in
H}c_\sigma\sigma(\mathfrak l),h_j(1+\tau)\sum_{\sigma\in
H}c_\sigma\sigma(\mathfrak l)\rangle.
\label{3.9}
\end{equation}
Calculating the right hand side of~\eqref{3.9} and taking into account that
\[
 \langle\sigma_1(\mathfrak l),\sigma_2(\mathfrak l)\rangle =
\begin{cases}
 1&\text{if $\sigma_1=\sigma_2,$}\\
0 & \text{if $\sigma_1\neq\sigma_2$}  
\end{cases}
\]
for any $\sigma_1,\sigma_2\in G$, we obtain
$\mathbf h^{-2}\{h_i(z'),h_j(z')\}= f_{ij}(c_1,\ldots,c_{n-1})$, where
$f_{ij}(c_1,\ldots,c_{n-1})$ is a quadratic form in $c_1,\ldots,c_{n-1}$ with
coefficients in $\mathbb Q$. Thus, $\det C$ is a polynomial $F$
in $c_1,\ldots,c_{n-1}$ of degree $2n$, and we wish to show that $F\neq0$ as a
polynomial. To do this we examine the coefficient at $c_1^{2n}$ in $F$.

Using bilinearness of the product $\langle\phantom x,\phantom x\rangle$, we obtain
from~\eqref{3.9} an equality $\mathbf h^{-2}\{h_i(z'),h_j(z')\}= A_{ij}+B_{ij}$, where
\[
 A_{ij}=\sum_{\beta,\gamma\in H}\langle h_ic_\beta\beta(\mathfrak
l),h_jc_\gamma\gamma(\mathfrak l)\rangle+\sum_{\beta,\gamma\in H}\langle h_i\tau
c_\beta\beta(\mathfrak l),h_j \tau c_\gamma\gamma(\mathfrak l)\rangle,
\]
\[
 B_{ij}=\sum_{\beta,\gamma\in H}\langle h_i\tau c_\beta\beta(\mathfrak l),h_j
c_\gamma\gamma(\mathfrak l)\rangle+\sum_{\beta,\gamma\in H}\langle
h_ic_\beta\beta(\mathfrak l),h_j\tau c_\gamma\gamma(\mathfrak l)\rangle.
\]
The elements $h_i\tau\beta$ and $h_j\gamma$ belong to the different classes modulo
$H$. The same statement holds for the elements $h_i\beta$ and $h_j\tau\gamma$.
Therefore, $B_{ij}=0$. Thus, $f_{ij}(c_1,\ldots,c_{n-1})=A_{ij}$.

To calculate $A_{ij}$, we note that
\begin{equation}
 \langle h_i c_\beta\beta(\mathfrak l), h_j c_\gamma\gamma(\mathfrak l)\rangle=
\begin{cases}
 c_\beta c_\gamma & \text{if $h_i\beta=h_j\gamma$,}\\
0 &\text{if $h_i\beta\neq h_j\gamma$.}
\end{cases}
\label{3.10}
\end{equation}
\begin{equation}
 \langle h_i\tau c_\beta\beta(\mathfrak l),h_j\tau c_\gamma\gamma(\mathfrak
l)\rangle=
\begin{cases}
 c_\beta c_\gamma &\text{if $h_i\tau\beta=h_j\tau\gamma$,}\\
0 & \text{if $h_i\tau\beta\neq h_j\tau\gamma$.} 
\end{cases}
\label{3.11}
\end{equation}
On the main diagonal of the matrix $C'$ man has $i=j$, that is $h_i=h_j$, so, it
follows from~\eqref{3.10} and~\eqref{3.11} that $A_{ii}=\sum_{k=1}^{n-1}2c_k^2$.
Thus, the product of all the elements on the main diagonal of the matrix $C'$
yields into $\det C'$ a summand
$\prod_{i=1}^nA_{ii}=(\sum_{k=1}^{n-1}2c_k^2)^n=2^nc_1^{2n}+\dots$, where dots
mean the sum of monomials, which contain $c_1$ in powers less than $2n$.

Now consider $A_{ij}$ for $i\neq j$, that is, in the case $h_i\neq h_j$. If this
is the case then the equalities $h_i\beta=h_j\gamma$ or
$h_i\tau\beta=h_j\tau\gamma$ can take place only for $\beta\neq \gamma$. Then it
follows from~\eqref{3.10} and~\eqref{3.11} that in the case $i\neq j$ the
quadratic form $f_{ij}(c_1,\ldots,c_{n-1})=A_{ij}$ does not contain $c_1^2$.
Therefore, any summand of $\det C'$, but the product of all the elements on the
main diagonal, contains $c_1$ in power at most $2n-1$. Thus, $\det
C'=F(c_1,\ldots, c_{n-1})$ is a non-zero polynomial in $c_1,\ldots,c_{n-1}$ with
coefficients in $\mathbb Q$. By Lemma~\ref{lemma3.1} the elements
$c_1,\ldots,c_{n-1}$ are algebraically independent over $\mathbb Q$. Therefore,
$\det C\neq 0$, and this means that $\mathfrak R_\ell(K)\neq0$. This concludes
the proof of the theorem.
\section{Proof of Theorem~3}
\setcounter{equation}{0}
Let $K$ be a finite Abelian extension of an imaginary quadratic $k=\mathbb
Q(\sqrt{-d})$ and $K$ is Galois
over $\mathbb Q$. Put $G=G(K/\mathbb Q),\quad H=G(K/k)$, and let $\tau$ be an
automorphism of complex conjugation.
 In particular, it means that $\tau\in G\setminus H$ and $\tau^2=1$.
Assume that prime $\ell$ splits completely in $K$, and let $S$ be the set of
all prime divisors of $\ell$ 
in $K$. Then $\ell$ splits in $k$ into a product of two prime divisors: 
$(\ell)={\mathfrak p}_1{\mathfrak p}_2$.
 So, man can present $S$ as a disjointed union $S=S_1\cup S_2$, where $S_i$ is
the set of all prime divisors 
of $K$ over $\mathfrak p_i,\quad i=1,2$. Thus, man obtains
 $D(K)=\mathbb Z_\ell[S]=D_1(K)\oplus D_2(K)$, where $D_i(K)=\mathbb
Z_\ell[S_i]$ for $i=1,2$. The group $H$
 acts on $S_1$ and $S_2$, and with respect to this action the modules $D_i(K)$
are isomorphic to
 $\mathbb Z_\ell[H]$ for $i=1,2$. The automorphism $\tau$ acts on $D(K)$,
interchanging $D_1(K)$ and $D_2(K)$.

Put $x,y\in D(K),\quad x=\sum_{i=1}^{2n}a_i\mathfrak l_i,\quad
y=\sum_{i=1}^{2n}b_i\mathfrak l_i$, where 
$n=[K:k],\quad \mathfrak l_1,\ldots,\mathfrak l_n\in S_1,\quad \mathfrak
l_{n+1},\ldots,\mathfrak l_{2n}\in 
S_2$. Then, but the pairing $\langle\phantom x,\phantom x\rangle$ that we have
deal with in the preceding
 section, we define two other pairings: $\langle
x,y\rangle_{K,1}=\sum_{i=1}^na_ib_i$ and
 $\langle x,y\rangle_{K,2}=\sum_{i=n+1}^{2n}a_ib_i$. We define also a relative
pairing
\begin{equation}
 \langle x,y\rangle_{K/k}=\langle x,y\rangle_{K,1}\mathfrak p_1+\langle
x,y\rangle_{K,2}
\mathfrak p_2\in D(k).
\label{4.0}
\end{equation}
An immediate checking shows that for any $\sigma\in G$ and any $x,y\in D(K)$ one
has
\begin{equation}
 \langle\sigma(x),\sigma(y)\rangle_{K/k}=\sigma(\langle x,y\rangle_{K/k}),
\label{4.1}
\end{equation}
where $G$ acts on $D(k)$ by restriction of automorphisms, that is $H$ acts on
$D(k)$ identically 
and $\tau$ interchanges $\mathfrak p_1$ and $\mathfrak p_2$. Obviously, for any
$x,y\in D(K)$
\begin{equation}
 \langle x,y\rangle_K=\langle x,y\rangle_{K,1}+\langle x,y\rangle_{K,2}.
\label{4.2}
\end{equation}
\begin{prop}
 Let $K,k,G,H,\tau,\mathfrak p_1,\mathfrak p_2$ be as above and $x\in D(K)$ such
an element that
 $\tau(x)=x$. Suppose that $\tau $ acts on $H$ by the rule $\tau(h)=\tau
h\tau^{-1}=\tau h\tau=h^{-1}$ for
 any $h\in H$. Then for any $h_i,h_j\in H$ an equality holds
\begin{equation}
 \tau\bigl(\langle h_i(x),h_j(x)\rangle_{K/k}\bigr)=\langle
h_i(x),h_j(x)\rangle_{K/k}.
\label{4.3}
\end{equation}
In other words,
\begin{equation}
 \langle h_i(x),h_j(x)\rangle_1=\langle h_i(x),h_j(x)\rangle_2.
\label{4.4}
\end{equation}
If $A$ is a $\mathbb Z_\ell[G]$-module generated by $x$ then for any $x,y\in A$
man has
 $\tau\bigl(\langle x,y\rangle_{K/k}\bigr)= \langle x,y\rangle_{K/k}$.
\label{prop4.1}
\end{prop}
{\bf Proof. } It follows from~\eqref{4.1} that
\begin{multline*}
 \tau\bigl(\langle h_i(x),h_j(x)\rangle_{K/k}\bigr)=\langle \tau h_i(x),\tau
h_j(x)\rangle_{K/k}= \\
=\langle \tau h_i\tau\tau(x),\tau h_j\tau\tau(x)\rangle_{K/k}=\langle
h_i^{-1}(x),h_j^{-1}(x)\rangle_{K/k}.
\end{multline*}
Applying an automorphism $h_i h_j$ and using~\eqref{4.1} once more, we get
\begin{multline*}
 \langle h_i^{-1}(x),h_j^{-1}(x)\rangle_{K/k}=h_ih_j\langle
h_i^{-1}(x),h_j^{-1}(x)\rangle_{K/k}= \\ 
=\langle h_j(x), h_i(x)\rangle_{K/k}=\langle h_i(x),h_j(x)\rangle_{K/k}.
\end{multline*}
This proves the proposition.

Let $H$ be a finite Abelian group and $\hat H$ the group of one-dimensional
characters of $H$ with values in $\bar{\mathbb Q}_\ell$. Let $\Phi$ be the set
of all irreducible characters of $H$ defined over $\mathbb Q_\ell$. Any
character $\varphi\in \Phi$ is a sum of all $\chi\in \hat H$, which are
conjugate by the action of the Galois Group $G(\bar{\mathbb Q}_\ell/\mathbb
Q_\ell)$. If $\chi$ enters such sum we say that $\chi$ enters $\varphi$ or
$\chi$ divides $\varphi$ and write down this fact by the formula
$\varphi=\sum_{\chi|\varphi}\chi$. If $e_\chi=|H|^{-1}\sum_{h\in H}\bar\chi(h)h$
is an idempotent that corresponds $\chi$ then the idempotent
$e_\varphi=\sum_{\chi|\varphi}e_\chi\in \mathbb Q_\ell[H]$ corresponds
$\varphi$.

Let $A\cong\mathbb Q_\ell[H]$ be a free $\mathbb Q_\ell[H]$-module of rank~1.
Then $A\cong \oplus_{\varphi\in \Phi}A_\varphi$, where $A_\varphi=e_\varphi A$.
For any $\varphi\in \Phi\quad A_\varphi\neq 0$ and $A_\varphi $ is an
irreducible $\mathbb Q_\ell[H]$-module. Note that for $\varphi\neq\varphi_1$ one
has $A_\varphi\ncong A_{\varphi_1}$.
For a character $\varphi\in \Phi$, where $\varphi=\sum_{\chi|\varphi}\chi$, by
$\bar\varphi$ we denote the character $\sum_{\chi|\varphi}\chi^{-1}$. Obviously,
if $\varphi$ is irreducible then so is $\bar\varphi$ and vice versa.
\begin{prop}
Let $A$ be a free $\mathbb Q_\ell[H]$-module of rank 1, which carries a
non-degenerate symmetric bilinear pairing $(\phantom x,\phantom x)\colon A\times
A\to\mathbb Q_\ell$, which satisfies the condition $(ha,hb)=(a,b)$ for any
$a,b\in A,\quad h\in H$. Then the pairing $(\phantom x,\phantom x)\colon
A_\varphi\times A_{\bar\varphi}\to\mathbb Q_\ell$ is non-degenerate, while the
pairing $(\phantom x,\phantom x)\colon A_\varphi\times A_\psi\to\mathbb Q_\ell$
vanishes for any $\varphi\in \Phi$ and $\psi\neq\bar\varphi$.
 \label{prop4.2}
\end{prop}
{\bf Proof. } Let $m$ be an exponent of $H$ and $L=\mathbb Q_\ell(\zeta_m)$,
where $\zeta_m$ is a primitive root of unity of degree $m$. Then any character
$\chi\in\hat H$ is defined over $L$ and the $L[H]$-module $B=A\otimes_{\mathbb
Q_\ell} L$ decomposes over $L$ into one-dimensional components
$B=\oplus_{\chi}B_\chi$, where $B_\chi=e_\chi B$. Here
$A_\varphi\otimes_{\mathbb Q_\ell}L=\oplus_{\chi|\varphi}B_\chi$ for any
$\varphi\in \Phi$. The pairing $(\phantom x,\phantom x)$ induces $L$-bilinear
pairing $(\phantom x,\phantom x)_1\colon B\times B\to L$, which again satisfies
condition $(a,b)_1=(ha,hb)_1$ for any $a,b\in B,\quad h\in H$. If $a,b\in A$
then $(a,b)=(a,b)_1$. Consider the pairing $(\phantom x,\phantom x)_1\colon
B_{\chi}\times B_{\theta}\to L$ for two characters $\chi,\theta\in \hat H$.
Since $he_\chi=\chi(h)e_\chi$ for any $h\in H,\quad \chi\in\hat H$, we get
\[
 (e_{\chi}a,e_{\theta}b)_1=(he_{\chi}a,he_{\theta}b)_1=
(\chi(h)e_{\chi}a,\theta(h)e_{\theta}b)_1=\chi\theta(h)(e_{\chi}a,e_{\theta}b)_1
\]
 for any $a,b \in B,\quad h\in H$.

Thus, if $\chi\theta\neq 1$ then $(B_{\chi},B_{\theta})_1=0$. If $\psi\neq
\bar\varphi$ then for any $\chi|\varphi,\theta|\psi$ one has $\chi\theta\neq1$,
that is, the components $B_\varphi$ and $B_\psi$ are mutually orthogonal.
Therefore $A_\varphi,A_\psi$ are also orthogonal with respect to the pairing
$(\phantom x,\phantom x)$. The pairing $A_\varphi\times A_{\bar\varphi}$ must be
non-degenerate, since otherwise it should be zero because of irreducibility of
$A_\varphi$. In such a case the pairing $A_\varphi\times A$ should be zero also,
but this contradicts the conditions of the proposition. This proves the
proposition.
\begin{theorem}
Let $K$ be an Abelian extension of an imaginary quadratic field $k=\mathbb
Q(\sqrt{-d})$ with Galois group  $H=G(K/k)$. Suppose that $K$ is Galois over
$\mathbb Q$ and prime $\ell$ splits completely in $K$. Let the complex
conjugation $\tau\in G=G(K/\mathbb Q)$ acts on $H$ by the rule $\tau(h)=\tau
h\tau^{-1}=h^{-1}$ for any $h\in H$. Let the exponent $m$ of $H$ be such that
the primitive root of unity $\zeta_m\in\bar{\mathbb Q}_\ell$ of degree $m$ is
conjugated with  $\zeta_m^{-1}$ by the action of the Galois group
$G(\bar{\mathbb Q}_\ell/\mathbb Q_\ell)$. Then we have $\mathfrak R_\ell(K)\neq
0$.
\end{theorem}
{\bf Proof. } By the consequence of Theorem~1 there is an element $\omega\in
U_{S,1}(K)$ such that $\tau(\omega)=\omega$ and the elements
$\{h(\omega)\}_{h\in H}$ generate in $U_{S,1}(K)$ a $\mathbb  Z_\ell$-submodule
of finite index, which we denote by $P$. By virtue of Theorem~1 $P$ is of
$\mathbb Z_\ell$-rank $r+1=|H|=n$, where $r$ is the rank of the group of units
$U(K)$ of $K$. Therefore $P\cong\mathbb Z_\ell[H]$ as $\mathbb Z_\ell[H]$-module.
Note that any $g\in G\setminus H$ is of the form $g=h\tau$ for some $h\in H$
hence $g(\omega)=h(\omega)$, that is, $P$ is also a $G$-module.

Since the Leopoldt conjecture is valid in $K$, by Prop.~\ref{prop3.1} the natural map
$\di\colon P\to D(K)$ is an injection. This map induces two homomorphisms:
$\di_1\colon P\to D_1(K)$ and $\di_2\colon P\to D_2(K)$, so that
$\di=\di_1\oplus\di_2$.

For technical convenience we pass from $\mathbb Z_\ell$-modules to $\mathbb
Q_\ell$-spaces. Put $\hat P=P\otimes_{\mathbb Z_\ell}\mathbb Q_\ell$. We define
$\hat D(K)$, $\hat D_1(K)$ and $\hat D_2(K)$ in analogous way. Then the
above-mentioned maps $\di,\di_1$ and $\di_2$ induce the maps, which we denote by
the same symbols, $\di\colon\hat P\to\hat D(K),\quad \di_1\colon\hat P\to\hat
D_1(K),\quad\di_2\colon \hat P\to\hat D_2(K)$. The pairing $\langle\phantom
x,\phantom x\rangle_K$ extends uniquely to the pairing $\hat D(K)\times \hat
D(K)\to\mathbb Q_\ell$, which we also denote by $\langle\phantom x,\phantom
x\rangle_K$. Analogously, we have the pairings
$\langle\phantom x,\phantom x\rangle_{K/k}\colon\hat D(K)\times\hat D(K)\to\hat
D(k)$ and
\[
 \langle\phantom x,\phantom x\rangle_{K,1}\colon\hat D_1(K)\times\hat D_1(K)\to
\mathbb Q_\ell,\quad
\langle\phantom x,\phantom x\rangle_{K,2}\colon\hat D_2(K)\times\hat D_2(K)\to
\mathbb Q_\ell.
\]
To prove the theorem, we have to check that the pairing $\langle\phantom
x,\phantom x\rangle_K$, being restricted to $\di(\hat P)$, induces a
non-degenerate pairing
\begin{equation}
\langle\phantom x,\phantom x\rangle_K\colon\di(\hat P)\times\di(\hat
P)\to\mathbb Q_\ell.
\label{4.6} 
\end{equation}
Let $\Phi$ be the set of all $\mathbb Q_\ell$-irreducible characters of $H$.
Suppose that the pairing~\eqref{4.6} has a non-trivial left kernel, for example. Then
this kernel has to contain the component $\di(\hat P_\varphi)$ for some
$\varphi\in \Phi$, that is, man has $\di(\hat P_\varphi)\times\di(\hat P)=0$. On
the other hand, since $\di$ is an injection, at least one of two conditions
$\di_1(\hat P_\varphi)\neq 0$ or $\di_2(\hat P_\varphi)\neq0$ holds. Suppose
$\di_1(\hat P_\varphi)\neq 0$. Since both $\hat D_1(K)$ and $\hat P$ are
isomorphic to $\mathbb Q_\ell[H]$, we obtain $\di_1(\hat P_\varphi)=\hat
D_1(K)_\varphi$. By the conditions of the theorem $\varphi=\bar\varphi$, hence
by Proposition~\ref{prop4.2} the non-degenerate pairing $\langle\phantom
x,\phantom x\rangle_{K,1}$ induces a non-degenerate pairing $\hat
D_1(K)_\varphi\times\hat D_1(K)_\varphi\to\mathbb Q_\ell$. Therefore the pairing
\begin{equation}
 \langle\phantom x,\phantom x\rangle_{K,1}\colon\di(\hat
P_\varphi)\times\di(\hat P_\varphi)\to \mathbb Q_\ell
\label{4.7}
\end{equation}
is also non-degenerate.

By the consequence of Theorem~1 the module $P$ satisfies the assertion of
Proposition~\ref{prop4.1}, that is, for any $a,b\in P$ we have $\langle
a,b\rangle_{K,1}=\langle a,b\rangle_{K,2}$. Thus, for any $a,b\in \hat
P_\varphi$, one has the equality $\langle \di(a),\di(b)\rangle_K=2\langle
\di(a),\di(b)\rangle_{K,1}$. So the non-degeneracy of~\eqref{4.7} yields the
non-degeneracy of the pairing $\langle\phantom x,\phantom
x\rangle_K\colon\di(\hat P_\varphi)\times\di(\hat P_\varphi)\to \mathbb Q_\ell$,
which contradicts the assumption that $\di(\hat P_\varphi)$ enters the kernel
of~\eqref{4.6}. The case $\di_2(\hat P_\varphi)\neq0$ can be treated
analogously. This proves the theorem.

We assume in Theorem~3  a primitive root of unity $\zeta_m$ of degree $m$ to be
conjugated with $\zeta_m^{-1}$ in the field $\bar{\mathbb Q}_\ell$, where $m$ is
the exponent of $H$. Obviously, the validity of this condition depends only on
$\ell$ and $m$. The next proposition characterize all pairs $(\ell,m)$, which
have this property.
\begin{prop}
 Let $\ell$ be a prime number and $m$ an arbitrary natural number. If
$\ell\neq2$ then we write $m$ in the form $m=2^a\ell^bc$ (if $\ell=2$ in the
form $m=2^ac$, where $c$ is prime to $2\ell$. Let $c_1$ be the product of all
different $p$ that divide $c$.

The root $\zeta_m\in\bar{\mathbb Q}_\ell$ is conjugated with $\zeta_m^{-1}$ over
$\mathbb Q_\ell$ if and only if one of the following conditions holds:
\begin{itemize}
\item[{\rm({\it i})}] $\ell\neq 2,\quad a=0,1$ and there is a natural number $r$
such that $\ell^r+1\equiv 0\pmod{c_1}$;
\item[{\rm({\it ii})}] $\ell\equiv 3\pmod 4,\quad a\geqslant 2$ and there is an
odd natural number $r$ such that $\ell^r+1\equiv0\pmod{2^ac_1}$;
\item[{\rm({\it iii})}] $\ell=2$ and there is a natural number $r$ such that
$\ell^r+1\equiv0\pmod{c_1}$.
\end{itemize}
\label{prop4.3}
\end{prop}
 {\bf Proof.} Suppose that there is $\sigma\in G(\bar{\mathbb Q}_\ell/\mathbb
Q_\ell)$ such that for given $m$ one has $\sigma(\zeta_m)=\zeta_m^{-1}$. For
$m=1,2$ we can assume $\sigma$ to be an identity automorphism. For $m>2$, the
condition $\sigma(\zeta_m)=\zeta_m^{-1}$ means that $\sigma$ is  an automorphism of order two
of the field $F=\mathbb Q_\ell(\zeta_m)$.

At first, suppose that $\ell\neq2$. The field $F$ is a free composite over
$\mathbb Q_\ell$ of $F_1=\mathbb Q_\ell(\zeta_{\ell^b})$ and $F_2=\mathbb
Q_\ell(\zeta_{2^ac})$. The extension $F_1/\mathbb Q_\ell$ is purely ramified and
has a Galois group isomorphic to $(\mathbb Z_\ell/\ell^b\mathbb Z_\ell)^\times$.
The automorphism $\sigma_1$, which corresponds to -1 under the last isomorphism,
sends $\zeta_{\ell^b}$ into $\zeta_{\ell^b}^{-1}$ and fixes any root of unity of
degree prime to $\ell$. The extension $F_2/\mathbb Q_\ell$ is unramified,
hence have a cyclic Galois group. Suppose that in $G(F_2/\mathbb Q_\ell)$ there
is an automorphism $\sigma_2$ such that
$\sigma_2(\zeta_{2^ac})=\zeta_{2^ac}^{-1}$. Then the automorphism
$\sigma=\sigma_1\sigma_2$ satisfies the condition
$\sigma(\zeta_m)=\zeta_m^{-1}$. Thus, if $\ell\neq2$ then  $\zeta_m$ and
$\zeta_m^{-1}$ are conjugated if and only if the roots $\zeta_{2^ac}$ and
$\zeta_{2^ac}^{-1}$ are conjugated in $F_2$.

1. If $a=0,1$ then $\zeta_{2^ac}$ and $\zeta_{2^ac}^{-1}$ are conjugated if and
only if $\zeta_c$ and $\zeta_c^{-1}$ are conjugated by some $\sigma_2$. Let $p$
be a prime divisor of $c$. Let $E=F_2^{\sigma_2}$ and $\bar E,\bar F_2$ be
residue fields of $E$ and $F_2$ respectively. If $|\bar E|=\ell^r$ then $|\bar
E^\times|=\ell^r-1$ and $|\bar F_2^\times|=\ell^{2r}-1=(\ell^r-1)(\ell^r+1)$. If
$\zeta_p\in F_2$ then $|\bar F_2^\times|\equiv0\pmod p$ but man has
$\sigma_2(\zeta_p)=\zeta_p^{-1}$. It means that $\zeta_p\notin E$, that is,
$\ell^r-1\not\equiv 0\pmod p$. Since $p$ is an arbitrary prime divisor of $c$,
we obtain $\ell^r+1\equiv0\pmod{c_1}$.

On the contrary, assume that $\ell^r+1\equiv0\pmod{c_1}$ and $a=0,1$. Take for
$E$ the unramified extension of $\mathbb Q_\ell$ of degree $r$, and let $E_2$ be
the unique quadratic unramified extension of $E$. The condition
$\ell^r+1\equiv0\pmod{c_1}$ implies that $\ell^r-1$ is prime to $c_1$, that is,
$\zeta_{c_1}\in E_2$ and $\zeta_p\notin E$ for any $p|c$. Since $\zeta_p\in
E_2$, man can obtain the field $F_2=E_2(\zeta_c)$ as follows. We already have
$\zeta_{c_1}\in E_2$. If $c_1=c$ then man can put $F_2=E_2$. If $c_1\neq c$ then
there is $p|c$ such that $c_2=pc_1|c$. Put $E_2'=E_2(\sqrt[p]{\zeta_{c_1}})$.
Then either $E_2'=E_2$ or $E_2'/E_2$ is a Kummer extension and hence has degree
$p$. If $c_2\neq c$ then $c$ has a divisor $c_3=qc_2$, where $q$ is a prime
divisor of $c$ (possibly, $q=p$). Then either the field
$E_2'':=E_2'(\sqrt[q]{\zeta_{c_2}})$ coincides with $E_2'$ or it is an extension of
degree $q$ of the latter.

Continuing in such a way, we get some extension $L$ of $E_2$ such that $L$
contains $\zeta_c$ and $L$ has odd degree over $E_2$. Therefore $[L:E]=2f$ for
some odd $f$, and there is a subfield $L'\subset L$ such that $L$ is quadratic
over $L'$. In particular, this means that $L$ is an extension of $E_2$ of some
odd degree $f$. Let $\sigma_2$ be the unique non-trivial automorphism of $L/L'$.
Since $L'$ has odd degree over $E$ and for any $p|\ell$ the primitive root of
unity $\zeta_p$ belongs to the quadratic extension $E_2$ of $E$, we obtain that
$L'$ has no roots $\zeta_p$ for $p|c$. It means that $N_{L/L'}(\zeta_c)=1$ or,
in other words $\sigma_2(\zeta_c)=\zeta_c^{-1}$. Then, putting $F_2=L$ and
$E=L'$, we get the desired.

2. Now we assume that $\ell\neq2$ and $a\geqslant2$. If $\ell\equiv1\pmod 4$
then the field $\mathbb Q_\ell$ contains a primitive root of unity $\zeta_4$ of
degree 4. If $\zeta_m$ would be conjugated with $\zeta_m^{-1}$ then $\zeta_4$
should be conjugated with $\zeta_4^{-1}$. But it is impossible because of the
condition $\zeta_4\in\mathbb Q_\ell$. Thus in the case $a\geqslant 2$ the roots
$\zeta_m$ and $\zeta_m^{-1}$ can be conjugated only if $\ell\equiv3\pmod 4$.
Moreover, if $\sigma_2$ is an automorphism of order two of  $F_2=\mathbb
Q_\ell(\zeta_{2^ac})$ then $\sigma_2(\zeta_4)=\zeta_4^{-1}$, that is,
$\zeta_4\notin E=F_2^{\sigma_2}$. It means that $\ell^r-1\equiv2\pmod 4$, where
$r=[E:\mathbb Q_\ell]$, that is, $E$ is an extension of $\mathbb Q_\ell$ of odd
degree. Just as in the paragraph~1, we obtain the congruence
$\ell^r+1\equiv0\pmod{c_1}$.

Let $2^s$ be the maximal degree of 2-power that contains in $F_2$. Note that
$s$ is the 2-adic exponent of $\ell^{2r}-1$. As it was shown above, the
2-component of the group of the roots of unity in $E$ consists of $\pm1$. On the
other hand, since the extension $F_2/E$ is unramified, the norm map $N\colon
F_2^\times\to E^\times$ induces a surjection of the multiplicative groups of the
residue fields. It means that $N(\zeta_{2^s})=-1$, that is,
$\sigma_2(\zeta_{2^s})=-\zeta_{2^s}^{-1}$. Thus the condition
$\sigma_2(\zeta_{2^a})=\zeta_{2^a}^{-1}$ holds only if $s>a$. This inequality is
equivalent to the condition  $\ell^r+1\equiv0\pmod{2^a}$. This proves the
necessity of the condition formulated in paragraph~2. Man can prove its
sufficiency, reasoning as in paragraph~1.

3. Now suppose $\ell=2$. In this case $m=2^ac$, where $(2,c)=1$. We have to
construct an unramified extension $F_2=\mathbb Q_2(\zeta_c)/\mathbb Q_2$ and an
automorphism of order two $\sigma_2$ of the field $F_2$ such that
$\sigma_2(\zeta_c)=\zeta_c^{-1}$. The same reasoning, as we used in paragraph~1,
shows that such $\sigma_2$ exists if and only if $\ell^r+1\equiv0\pmod{c_1}$ for
some $r$. This proves the proposition.

{\bf Examples.} 1. If $\ell=5$ then $\ell^3+1=126=2\cdot 3^2\cdot 7$, hence any
$m$ of the form $m=2^a3^{s_1}5^{s_2}7^{s_3}$, where $a=0,1$, satisfies
condition~(i) of Prop.~\ref{prop4.3}. Therefore, for $\ell=5$ and any $K$ such
that the Galois group $H=G(K/k)$ is of exponent $m$ of the above mentioned form,
man has $\mathfrak R_5(K)\neq0$ provided 5 splits completely in $K$. 

2. Put $\ell=11$. Then $\zeta_5\in\mathbb Q_{11}$ hence Theorem~3 does not hold
for any $m$ that divides by~5.

Now we assume that the field $K$ from Theorem~3 is fixed. Our goal is to
characterize all prime $\ell$ such that Theorem~3 holds for the pair $(K,\ell)$.
\begin{prop}
 Let a field  $K$ be Galois over $\mathbb Q$ and $K$ is an Abelian extension of
some imaginary quadratic field $k$, as in Theorem~3. Let $m$ be an exponent of
$H=G(K/k)$ and $T$ the set of all primes that do not divide $m$. Put $F=\mathbb
Q(\zeta_m)$. 

If $E:=K\cap F$ is totally real (in particular, if $E=\mathbb Q$) then there are
infinitely many primes $\ell\in T$, which satisfy all conditions of Theorem~3.
In particular, man has $\mathfrak R_\ell(K)\neq0$ for these $\ell$. Otherwise,
the conditions of Theorem~3 take place for no $\ell\in T$.
\label{prop4.4}
\end{prop}
{\bf Proof. } Put $L=K\cdot F$. Let $\delta$ be the automorphism of complex
conjugation in $F$. If the field $E$ is totally real then $\delta$ acts on $E$
as identity. Since $L$ is a free composite of $K$ and $F$ over $E$, there is an
automorphism $\sigma\in G(L/E)$ such that $\sigma|_K={\rm id}$ and
$\sigma|_F=\delta$. By the Chebotarev density theorem there are infinitely many
prime divisors $\mathfrak l$ of $L$ such that $\mathfrak l$ is unramified in $L$
and the Frobenius automorphism $\varphi(\mathfrak l)$ equals $\sigma$. Let
$\ell$ be a prime rational number under $\mathfrak l$. Since $\sigma$ induces
the identity automorphism of $K$, we see that $\ell$ splits completely in
$K$; while the fact that $\sigma$ and $\delta$ coincide, being restricted to the
cyclotomic field $F=\mathbb Q(\zeta_m)$, means that $\ell\equiv-1\pmod{m}$. By
Prop.~\ref{prop4.3} any such $\ell$ satisfies conditions of Theorem~3 (with
$r=1$), hence for such $\ell$ one has $\mathfrak R_\ell(K)\neq0$.

On the contrary, assume that $E$ is an imaginary field. Suppose that for some prime $\ell\in T$ 
the roots $\zeta_m$ and $\zeta_m^{-1}$ are conjugated over $\mathbb
Q_\ell$. Then by Prop.~\ref{prop4.3} man has a congruence
$\ell^r+1\equiv0\pmod m$ for some $r$. Let $\varphi(\mathfrak l)$ be the Frobenius automorphism
corresponding to $\ell$. Then the last congruence means that
$\varphi(\mathfrak l)^r$ coincides with $\delta$ on $F$. Since $E$ is imaginary,
$\delta$ acts non-trivially on $E$, hence $\varphi(\mathfrak
l)$ acts non-trivially on $K$, that is, $\ell$ does not split completely in
$K$. This proves the proposition.

{\bf Remark.} Let $p\equiv3\pmod 4$ be a prime number, $k=\mathbb
Q(\sqrt{-p})$ and $K$ an Abelian extension of $k$, as in Theorem~ 3, such that
the group $H=G(K/k)$ is of exponent $p$. In this case we have $F=\mathbb Q(\zeta_p)$ and
$F\cap K\supset k$, so by Prop.~\ref{prop4.4} the assumptions of Theorem~3 
don't hold for any prime $\ell$ but $p$.

Nevertheless if $\ell=p$ one can apply this theorem and obtain $\mathfrak
R_\ell(K)\neq0$.

L.~V.~Kuz'min\\
NRC ``Kurchatovskii institute''\\
E-mail: helltiapa@mail.ru

\end{document}